\documentclass[%
 preprint,
 amsmath,amssymb,
 aps,
]{revtex4-2}

\usepackage{rotating}
\usepackage{makecell}
\usepackage{graphicx}
\usepackage{dcolumn}
\usepackage{bm}
\usepackage{amsthm} 
\theoremstyle{remark}
\newtheorem{remark}{Remark}

\begin{document}


\title{A novel time integration scheme for linear parabolic PDEs}

\author{Subhankar Nandi}
\author{Satyajit Pramanik}%
 \email{satyajitp@iitg.ac.in}
\affiliation{%
 Department of Mathematics, Indian Institute of Technology Guwahati, Guwahati - 781039, Assam, India
}%




\date{\today}

\begin{abstract}
This paper presents a class of Crank–Nicolson (CN) type schemes enhanced by radial basis function (RBF) interpolation for the time integration of linear parabolic partial differential equations (PDEs). The resulting RBF–CN schemes preserve the structural simplicity of the classical CN method while enabling higher-order temporal accuracy through the optimization of the shape parameter. Consistency analysis shows that the schemes are always at least second-order accurate and a simple choice of the shape parameter increases the accuracy by two orders over the standard CN scheme. Further optimization can reduce the next leading error term to attain even higher-order accuracy. A von Neumann stability analysis confirms that the stability conditions are essentially the same as those of the standard CN scheme. Several numerical experiments in 1D are carried out to verify the theoretical results under different boundary conditions. Particular focus is given to the startup stage, where several strategies for computing the initial steps are examined, and Gauss–Legendre implicit Runge–Kutta (IRK) methods are found to be the most effective. The experiments further demonstrate that, with optimal initialization, the proposed schemes deliver accuracy comparable to implicit Runge–Kutta methods while achieving nearly two orders of magnitude reduction in computational cost.
\end{abstract}

\maketitle


\section{Introduction} \label{intro}

Accurate and efficient time integration schemes are essential in the numerical simulation of time-dependent partial differential equations (PDEs), which model a wide array of physical phenomena across science and engineering applications. Over the years, many \citep{curtiss1952integration, calahan1968stable,boscarino2016high, sebastiano2023high} time-stepping methods have been developed, each balancing computational cost, stability, and accuracy in different ways. These methods are generally categorized as explicit, implicit, or semi-implicit \citep{butcher2016numerical}. Explicit schemes are easy to implement and computationally inexpensive per time step, yet they require very small step sizes to maintain numerical stability. However, their stability is often restricted by stringent time step constraints, especially for stiff systems. On the other side, the implicit schemes allow for larger time steps without violating stability constraints, though at the expense of increased computational complexity due to the need to solve systems of equations at each step. 

Among the widely used approaches, certain schemes are valued for their ability to remain stable regardless of the chosen time step size, particularly useful in stiff problems or long-time simulations. One such classical choice is the standard Crank–Nicolson (CN) time integration scheme, which has gained widespread use due to its simplicity, unconditional stability, and second-order temporal accuracy. Like many classical schemes, CN and related methods are typically derived through polynomial interpolation or Taylor expansions, where the attainable order of accuracy is inherently tied to the degree of the underlying polynomial approximation. While higher-order extensions of these polynomial-based methods are theoretically possible, in practice they often compromise stability, increase computational cost, or diminish the simplicity that makes lower-order schemes attractive. Such challenges have motivated the search for alternative time integration frameworks that can systematically achieve higher-order accuracy without sacrificing stability or increasing algorithmic complexity.

A common approach to improving the local accuracy of a scheme is to incorporate the information about the solution derivatives estimated to a certain order. However, the schemes based on polynomial interpolation offer limited flexibility in this regard. In contrast, radial basis function (RBF) interpolation \citep{buhmann2000radial} provides a flexible and adaptive framework that is capable of delivering arbitrarily high-order accuracy. This property motivates a modification in the CN scheme using RBF interpolation in time. A similar idea has been used in many studies to construct higher-order versions of classical schemes such as Euler’s method, midpoint methods, explicit Adams-Bashforth methods \citep{gu2020adaptive, gu2021adaptive}, and explicit Runge-Kutta (RK) methods of various orders \citep{gu2024explicit} for solving initial value problems (IVPs). This approach was also applied to modify the ENO and WENO methods \citep{gu2017radial} for solving hyperbolic PDEs. The common feature in all these works is the use of infinitely smooth RBFs for interpolation. These RBFs include a shape parameter, which plays a key role in the modification of the methods. In these studies, an equation involving the shape parameter $\epsilon$ is derived from the local truncation error. It is then suggested that solving this equation gives the optimal value of $\epsilon$, which improves the overall accuracy when used in the modified scheme.

However, most of these works mainly focus on the modification of explicit methods and show that, with a simple choice of the shape parameter $\epsilon$, which usually comes from the leading term of the local truncation error, the modified scheme becomes one order more accurate than the classical version. Although it is possible to choose $\epsilon$ in other ways to obtain schemes of even higher (arbitrary) order, this increases the complexity of the method and makes it less effective for practical use.

In this study, we develop a higher-order version of the implicit time integration scheme like CN using a similar strategy for linear parabolic PDEs. For the interpolation, we choose three infinitely smooth RBFs, namely Gaussian, multi-quadric (MQ), and inverse multi-quadric (IMQ). To construct the scheme, we first consider two consecutive time nodes and construct an interpolating function using RBF interpolation, where the values of the solution at these two nodes are used as constraints. We then compute the time derivative of the interpolant at both nodes. By using the derivative at the left node in the PDE, we obtain a modified version of the forward Euler method \citep{gu2020adaptive, gu2021adaptive}. Similarly, using the derivative at the right node gives us a modified version of the backward Euler method. Taking the average of these two modified schemes, we found a modified form of the CN scheme, which we call as RBF-CN scheme. One advantage of these schemes is that they retain the structure of the original CN scheme and recover it precisely in the limit when the shape parameter vanishes.

We next examine the consistency of the developed RBF-CN schemes in time by assuming the spatial operator is exact. We see that the leading-order in the local truncation error is two, which represents the actual order of the scheme in time. Thus, the RBF-CN schemes are theoretically at least second-order accurate for any choice of $\epsilon$. For higher-order accuracy, we locally optimize the shape parameter so that the coefficients of the leading error terms vanish. We derive the expressions for obtaining optimal $\epsilon$ for different orders of convergence in terms of the solution variable. We found that the choice of optimal $\epsilon$ that vanishes the first leading order error terms, makes the scheme two orders more accurate than the standard CN. In addition to consistency, we analyze the stability of the proposed schemes using the von Neumann approach for a general linear operator and discuss specific cases. We find that the schemes exhibit stability conditions similar to those of the classical CN method applied to a general operator.

To verify the theoretical findings, we test the proposed RBF–CN schemes on problems with different boundary conditions. In the experiments, we first examine the start-up difficulty that appears because initial step solutions are required to approximate the optimal $\epsilon$. For the case of fourth-order accuracy, we consider several approaches for computing the initial steps, including CN with the same step size, CN with a refined step size, CN with Richardson extrapolation, and the Gauss–Legendre implicit Runge–Kutta (IRK) method \citep{butcher1964implicit}. Among these, the IRK method gives the most accurate results while also being efficient in terms of cost. We then study the effect of different orders of approximation for the optimal $\epsilon$ and point out suitable choices. Finally, we compare the full IRK method with the RBF–CN scheme that uses IRK for the initial steps. The results show that the proposed scheme achieves similar accuracy but is nearly two orders of magnitude faster than IRK.

The structure of the paper is organized as follows. Section \ref{sec:RBFs} introduces the basics of RBF interpolation. Section \ref{sec:RBF-CN} presents the formulation of the RBF–CN schemes. Section \ref{sec:consistency} establishes their consistency properties, while \S \ref{sec:stability} addresses stability analysis. Section \ref{sec:numerical} demonstrates the accuracy and efficiency of the methods through numerical experiments before summarizing the concluding remarks in \S \ref{sec:conclusion}. 

\section{Interpolation using RBFs}\label{sec:RBFs}
We begin by introducing the general formulation of the RBF interpolation in one dimension. Let $u(t)$ be a sufficiently smooth function defined on an interval $[t_0, t_N]$, where $t_0, \dots, t_N$ are $N+1$ distinct nodes. To approximate $u(t)$ at any point in the interval based on its values at these nodes, we construct an interpolant of the form
\begin{equation}\label{RBF_interpolant}
    u(t) \approx \sum_{j=0}^N \lambda_j\, \phi_j(|t - t_j|) + p(t),
\end{equation}
where $\phi_j(t) = \phi(|t - t_j|)$ is the radial basis function with $t_j$ as the center, $p(t)$ is a low-order polynomial of degree at most $m-1$ with $m\le N$ included to ensure a consistent system and $\lambda_j$ are coefficients to be determined. Incorporating the interpolation constraints $u(t_i) = u_i$ for each $i = 0, \dots, N$, and additional side conditions to ensure a unique solution, a linear system for the coefficients is obtained. In this study, we omit the polynomial term (i.e., set $p(t) = 0$) for simplicity. The resulting linear system for the unknown coefficients $\lambda_j$ can then be written in matrix form as
\begin{equation}\label{matrixsystem}
A \boldsymbol{\lambda} = \boldsymbol{u}.
\end{equation}
where $A_{ij} = \phi(|t_i - t_j|)$, $\boldsymbol{\lambda} = [\lambda_0, \dots, \lambda_N]^T$ contains the RBF coefficients and $\mathbf{u} = [u_0, \dots, u_N]^T$ contains the known function values at the interpolation nodes. 

There are generally two categories of RBFs commonly used for interpolation: those without a shape parameter and those with one. When no shape parameter is involved, the interpolation system can be solved directly. In contrast, when a shape parameter $\epsilon$ is present, it must be prescribed either as a global constant or through optimization to improve interpolation accuracy. The value of the parameter $\epsilon$ may be real or complex, depending on the error optimization, and can also vary across all nodes. In this work, we proceed with infinitely smooth RBFs with a shape parameter that is identical across all interpolation nodes.

\subsection{Gaussian RBF interpolation}
We first consider the Gaussian RBF, defined as:
\begin{equation}\label{Gaussian_RBF}
    \phi_k(t) = e^{-\epsilon^2 (t - t_k)^2}.
\end{equation}
For interpolation using only two points, $(t_0, u_0)$ and $(t_1, u_1)$, the corresponding linear system \eqref{matrixsystem} becomes
\begin{equation}
\begin{bmatrix}
1 & e^{-\epsilon^2 \Delta t^2} \\
e^{-\epsilon^2 \Delta t^2} & 1
\end{bmatrix}
\begin{bmatrix}
\lambda_0 \\
\lambda_1
\end{bmatrix}
=
\begin{bmatrix}
u_0 \\
u_1
\end{bmatrix},
\end{equation}
where $\Delta t = t_1 - t_0$. The resulting interpolant is given by
\begin{equation}\label{Gaussian_RBF_interpolant}
u(t) = \lambda_0\, e^{-\epsilon^2 (t - t_0)^2} + \lambda_1\, e^{-\epsilon^2 (t - t_1)^2},
\end{equation}
with
\begin{equation}\label{Gaussian_RBF_coefficients}
\lambda_0 = \frac{u_0 - u_1 \,e^{-\epsilon^2 \Delta t^2} }{1 - e^{-2\epsilon^2 \Delta t^2}}, \qquad
\lambda_1 = \frac{u_1 - u_0 \, e^{-\epsilon^2 \Delta t^2} }{1 - e^{-2\epsilon^2 \Delta t^2}}.
\end{equation}

\subsection{MQ RBF interpolation}
For the interpolation using MQ RBF, given by  
\begin{equation}
    \phi_k(t) = \sqrt{1 + \epsilon^2 (t - t_k)^2},
\end{equation}  
the interpolant becomes
\begin{equation}
\label{MQ_RBF_interpolant}
    u(t) = \lambda_0\, \sqrt{1 + \epsilon^2 (t - t_0)^2} + \lambda_1\, \sqrt{1 + \epsilon^2 (t - t_1)^2},
\end{equation}
with
\begin{equation}\label{MQ_RBF_coefficients}
\lambda_0 = -\frac{u_0 - u_1 \, \sqrt{1 + \epsilon^2 \Delta t^2} }{ \epsilon^2 \Delta t^2}, \qquad
\lambda_1 = \frac{  u_0 \, \sqrt{1 + \epsilon^2 \Delta t^2} - u_1}{ \epsilon^2 \Delta t^2}.
\end{equation}

\subsection{IMQ RBF interpolation}
Similarly, for the IMQ RBF defined as
\begin{equation}
    \phi_k(t) = \frac{1}{\sqrt{1 + \epsilon^2 (t - t_k)^2}},
\end{equation}  
the resulting interpolant is
\begin{equation}\label{IMQ_RBF_interpolant}
u(t) = \lambda_0\, \frac{1}{\sqrt{1 + \epsilon^2 (t - t_0)^2}} + \lambda_1\, \frac{1}{\sqrt{1 + \epsilon^2 (t - t_1)^2}},
\end{equation}
with
\begin{equation}\label{IMQ_RBF_coefficients}
\lambda_0 = \frac{u_0 \, (1 + \epsilon^2 \Delta t^2) -  u_1 \, \sqrt{1 + \epsilon^2 \Delta t^2} }{\epsilon^2 \Delta t^2}, \qquad
\lambda_1 = \frac{u_1 \,(1 + \epsilon^2 \Delta t^2)  -   u_0 \,\sqrt{1 + \epsilon^2 \Delta t^2} }{\epsilon^2 \Delta t^2}.
\end{equation}

\section{Development of the RBF-CN Scheme}\label{sec:RBF-CN}
We consider the linear parabolic PDE of the form
\begin{equation}\label{eq:linearPDE}
\frac{\partial u(\mathbf{x}, t)}{\partial t} = \mathcal{L} u(\mathbf{x}, t),
\end{equation}
where $ \mathbf{x} \in \mathbb{R}^d $,  $t \geq 0 $, and  $\mathcal{L}$ is a linear spatial differential operator. We aim to develop an efficient, high-order accurate CN-type time-stepping scheme in the temporal direction. The key idea is to replace the classical time discretization of standard finite difference methods with one constructed from a temporal interpolant using RBF. The development of the scheme may require the solution to be sufficiently smooth, which will be discussed later on. However, we do not address the spatial discretization in the development of the scheme, except where it is necessary to ensure stability.

\subsection{Construction using Gaussian RBF}
We begin with the Gaussian RBF interpolant defined in \eqref{Gaussian_RBF_interpolant}. Differentiating with respect to $t$ and evaluating at $t = t_0$ yields
\begin{equation}
\left( \frac{du}{dt} \right)_{t = t_0} = \lambda_0 \phi_0'(t_0) + \lambda_1 \phi_1'(t_0),
\quad \text{where} \quad \phi_k'(t) = -2\epsilon^2(t - t_k) e^{-\epsilon^2(t - t_k)^2}.
\end{equation}
This simplifies to
\begin{equation}
\left( \frac{du}{dt} \right)_{t = t_0} = \frac{u_1 - u_0 \, e^{-\epsilon^2 \Delta t^2} }{1 -e^{-2\epsilon^2 \Delta t^2}} \cdot 2 \epsilon^2 \Delta t \, e^{-\epsilon^2 \Delta t^2}.
\end{equation}
Substituting into the PDE \eqref{eq:linearPDE} gives
\begin{equation}\label{eq:GS_Forward1}
u_1 - e^{-\epsilon^2 \Delta t^2} u_0 = \frac{1 - e^{-2\epsilon^2 \Delta t^2}} {2 \epsilon^2 \Delta t \, e^{-\epsilon^2 \Delta t^2}} \, \mathcal{L}u_0,
\end{equation}
which can be expressed equivalently as
\begin{equation}\label{eq:GS_Forward2}
u_1 = u_0 e^{-\epsilon^2 \Delta t^2} + \frac{e^{\epsilon^2 \Delta t^2} - e^{-\epsilon^2 \Delta t^2}} {2 \epsilon^2 \Delta t} \, \mathcal{L}u_0.
\end{equation}

Similarly, applying the interpolation in the interval $[t_1, t_2]$ using the points $(t_1, u_1)$ and $(t_2, u_2)$, and computing the first derivative at $t = t_1$ with the same step size $\Delta t = t_2 - t_1$, yields an analogous expression for $u_2$. Repeating the same argument over the interval $[t_n, t_{n+1}]$ and allowing the shape parameter to vary with $n$, yields
\begin{equation}\label{eq:GS_Forward}
u_{n+1} = u_n e^{-\epsilon_n^2 \Delta t^2} + \frac{e^{\epsilon_n^2 \Delta t^2} - e^{-\epsilon_n^2 \Delta t^2}} {2 \epsilon_n^2 \Delta t} \, \mathcal{L}u_n.
\end{equation}
We refer to equation \eqref{eq:GS_Forward} as the Gaussian RBF forward Euler scheme.
In the special case where the shape parameter tends to zero, i.e., $\epsilon \to 0$, 
the prefactor of $\mathcal{L}(u_n)$ in \eqref{eq:GS_Forward} satisfies
\begin{equation}
\lim_{\epsilon \to 0} 
\frac{e^{\epsilon^2 \Delta t^2} - e^{-\epsilon^2 \Delta t^2}}{2 \epsilon^2 \Delta t}
= \Delta t,
\end{equation}
and equation \eqref{eq:GS_Forward} reduces to the classical forward Euler scheme
\begin{equation}\label{eq:standard_forward}
u_{n+1} = u_n + \Delta t\, \mathcal{L}u_n.
\end{equation}

Applying again the same interpolation but evaluating at $t= t_1$ gives
\begin{equation}
    \left(\frac{du}{dt}\right)_{t = t_1} = \frac{u_1 e^{-\xi^2 \Delta t^2} - u_0}{1 - e^{-2\xi^2 \Delta t^2}} \cdot 2 \epsilon^2 \Delta t \, e^{-\epsilon^2 \Delta t^2}
\end{equation}
Substituting into the PDE \eqref{eq:linearPDE} yields 
\begin{equation}\label{eq:GS_Backward1}
u_1 \, e^{-\epsilon^2 \Delta t^2} -  u_0 = \frac{1 - e^{-2\epsilon^2 \Delta t^2}} {2 \epsilon^2 \Delta t \, e^{-\epsilon^2 \Delta t^2}} \, \mathcal{L}u_1,
\end{equation}
which simplifies to 
\begin{equation}\label{eq:GS_Backward2}
u_1  = u_0 \, e^{\epsilon^2 \Delta t^2} +  \frac{e^{2\epsilon^2 \Delta t^2} - 1} {2 \epsilon^2 \Delta t} \, \mathcal{L}u_1,
\end{equation}
The general $n$-step form is
\begin{equation}\label{eq:GS_Backward}
u_{n+1}  = u_n \, e^{\epsilon^2_n \Delta t^2} +  \frac{e^{2\epsilon^2_n \Delta t^2} - 1} {2 \epsilon^2_n \Delta t} \, \mathcal{L}u_{n+1},
\end{equation}
which we refer to as the Gaussian RBF backward Euler scheme. In the limit $\epsilon \to 0$, the coefficient multiplying $\mathcal{L}u_{n+1}$ in \eqref{eq:GS_Backward} approaches $\Delta t$, which recovers the original standard backward Euler scheme
\begin{equation}\label{eq:standard_backward}
u_{n+1} = u_n + \Delta t\, \mathcal{L}u_{n+1}.
\end{equation}
Note that this scheme \eqref{eq:GS_Backward} is implicit, while the Gaussian RBF forward scheme \eqref{eq:GS_Forward} was explicit as similar to the standard forward and backward Euler scheme.

Adding the forward \eqref{eq:GS_Forward} and backward expressions \eqref{eq:GS_Backward} at the same time level eliminates the exponential scaling on $u_0$ and $u_{1}$, giving
\begin{equation}
\big(u_1 - u_0\big)\left( 1 + e^{-\epsilon^{2} \Delta t^{2}} \right)
= \frac{1 - e^{-2 \epsilon^{2} \Delta t^{2}}}
       { 2 \epsilon^{2} \Delta t \, e^{-\epsilon^{2} \Delta t^{2}} }
  \left( \mathcal{L}u_0 + \mathcal{L}u_{1} \right)
\end{equation}
which can be rearranged to give
\begin{equation}\label{eq:gaussian_CN1}
u_1 = u_0 + \frac{e^{\epsilon^2 \Delta t^2} - 1}{2 \epsilon^2 \Delta t} \left( \mathcal{L}u_0 + \mathcal{L}u_{1} \right)
\end{equation}
Extending this formulation to all time levels yields
\begin{equation}\label{eq:gaussian_CN}
u_{n+1} = u_{n} + \frac{e^{\epsilon^2_n \Delta t^2} - 1}{2 \epsilon^2_n \Delta t} \left( \mathcal{L}u_n + \mathcal{L}u_{n+1} \right)
\end{equation}
which we refer to as the Gaussian RBF CN scheme. When the shape parameter tends to zero, $\epsilon \to 0$, the prefactor of $\mathcal{L}$ approaches to $\Delta t /2$ and recovers the standard CN scheme
\begin{equation}\label{eq:standard_CN}
    u_{n+1} = u_{n} + \frac{\Delta t}{2} \left( \mathcal{L}u_n + \mathcal{L}u_{n+1} \right).
\end{equation} 

\subsection{Construction using MQ RBF}

For the MQ RBF, the expression for the temporal derivative at $t=t_0$ simplifies to
\begin{equation}
\left( \frac{du}{dt} \right)_{t = t_0} = \frac{u_1 - m u_0}{m \, \Delta t},
\quad \text{where} \quad m = \sqrt{1 + \epsilon^2 \Delta t^2}.
\end{equation}
Substituting into \eqref{eq:linearPDE} gives the MQ analogue of the forward update:
\begin{equation}\label{eq:MQ_Forward1}
u_1 - m u_0 = m \Delta t \; \mathcal{L}u_0,
\end{equation}
or, equivalently,
\begin{equation}\label{eq:MQ_Forward2}
u_1 = m \left( \, u_0 + m \Delta t \; \mathcal{L}u_0 \right).
\end{equation}
Allowing $\epsilon$ (and hence $m$) to vary with the time step $n$ yields the general $n$-step form
\begin{equation}\label{eq:MQ_Forward}
u_{n+1} = m_n \left( \, u_n + m_n \Delta t \; \mathcal{L}u_n \right),
\qquad m_n = \sqrt{1 + \epsilon_n^2 \Delta t^2},
\end{equation}
which we refer to as the MQ RBF forward Euler scheme. In the limit $\epsilon \to 0$, we have $m \to 1$ and \eqref{eq:MQ_Forward} recovers the classical forward Euler scheme.

Repeating the interpolation but differentiating at $t=t_1$ leads to the backward form. Evaluating the derivative at $t = t_1$ gives
\begin{equation}
\left( \frac{du}{dt} \right)_{t = t_1}
= \frac{m u_1 - u_0}{m \, \Delta t},    
\end{equation}
and substitution into \eqref{eq:linearPDE} yields
\begin{equation}\label{eq:MQ_Backward1}
m u_1 - u_0 = m \Delta t \; \mathcal{L}u_1,
\end{equation}
which simplifies to
\begin{equation}\label{eq:MQ_Backward2}
u_1 = \frac{u_0}{m} + \Delta t \; \mathcal{L}u_1.
\end{equation}
The general implicit $n$-step form is
\begin{equation}\label{eq:MQ_Backward}
u_{n+1} = \frac{u_n}{m_n} + \Delta t \; \mathcal{L}u_{n+1},
\end{equation}
which we call the MQ RBF backward Euler scheme. As $\epsilon \to 0$, $m \to 1$ and \eqref{eq:MQ_Backward} reduces to the classical backward Euler method.

Adding \eqref{eq:MQ_Forward} and \eqref{eq:MQ_Backward} eliminates the multiplicative scalings on $u_n$ and $u_{n+1}$ and yields the MQ CN form:
\begin{equation}
\big(u_{n+1} - u_n\big) (1 + m_n)
= m_n \Delta t \left( \mathcal{L}u_n + \mathcal{L}u_{n+1} \right),
\end{equation}
which can be rearranged as
\begin{equation}\label{eq:MQ_CN}
u_{n+1} = u_n + \frac{m_n}{1 + m_n}\, \Delta t \left( \mathcal{L}u_n + \mathcal{L}u_{n+1} \right).
\end{equation}
Note that $\displaystyle \frac{m_n}{1+m_n}\to \tfrac12$ as $\epsilon_n\to 0$, so \eqref{eq:MQ_CN} recovers the standard CN scheme in the $\epsilon\to 0$ limit. Note that the MQ-based schemes are similar but algebraically simpler than the Gaussian case as exponentials reduce to the multiplier $\sqrt{1 + \epsilon^2 \Delta t^2}$. 

\subsection{Construction using IMQ RBF}

Consider IMQ RBF interpolant given in \eqref{IMQ_RBF_interpolant}. Evaluating the temporal derivative of the interpolant at $t=t_0$ gives
\begin{equation}
\left(\frac{du}{dt}\right)_{t=t_0} = \frac{u_1 - \tfrac{1}{m} u_0}{1 - \tfrac{1}{m^2}} \cdot \frac{\epsilon^2 \Delta t}{m^3}, \qquad \text{with} \quad  m = \sqrt{1 + \epsilon^2 \Delta t^2}.
\end{equation}
Since $1 - \tfrac{1}{m^2} = \frac{m^2 - 1}{m^2} = \frac{\epsilon^2 \Delta t^2}{m^2}$, this simplifies to
\begin{equation}
\left(\frac{du}{dt}\right)_{t=t_0}
= \frac{m\,(u_1 - \tfrac{1}{m} u_0)}{\Delta t}.
\end{equation}

Substituting into the PDE \eqref{eq:linearPDE} yields the IMQ forward scheme
\begin{equation}\label{eq:IMQ_Forward1}
m\,u_1 - u_0 = m\,\Delta t \; \mathcal{L}u_0,
\end{equation}
which can be written equivalently as
\begin{equation}\label{eq:IMQ_Forward2}
u_1 = \frac{u_0}{m} + \Delta t\; \mathcal{L}u_0.
\end{equation}
Allowing the shape parameter to vary with the time step gives the general $n$-step IMQ forward scheme:
\begin{equation}\label{eq:IMQ_Forward}
u_{n+1} = \frac{u_n}{m_n} + \Delta t\; \mathcal{L}u_n.
\end{equation}
We refer to \eqref{eq:IMQ_Forward} as the IMQ RBF forward Euler scheme.

Similarly, differentiating the interpolant at $t=t_1$ yields
\begin{equation}
\left(\frac{du}{dt}\right)_{t=t_1} = \frac{\tfrac{1}{m} u_1 - u_0}{\Delta t}\, m.
\end{equation}
Substituting into \eqref{eq:linearPDE} leads to
\begin{equation}\label{eq:IMQ_Backward1}
u_1 = m\,u_0 + \Delta t\; \mathcal{L}u_1.
\end{equation}
The general implicit form is
\begin{equation}\label{eq:IMQ_Backward}
u_{n+1} = m_n\,u_n + \Delta t\; \mathcal{L}u_{n+1},
\end{equation}
which we call the IMQ RBF backward Euler scheme. 

Adding \eqref{eq:IMQ_Forward1} and \eqref{eq:IMQ_Backward1} eliminates the asymmetric coefficients on $u_n$ and $u_{n+1}$ and yields the IMQ Crank--Nicolson form:
\begin{equation}
\big(u_{n+1} - u_n\big)(1 + m_n) = \Delta t \left( \mathcal{L}u_n + \mathcal{L}u_{n+1} \right),
\end{equation}
which can be rearranged to
\begin{equation}\label{eq:IMQ_CN}
u_{n+1} = u_n + \frac{\Delta t}{1+m_n} \left( \mathcal{L}u_n + \mathcal{L}u_{n+1} \right).
\end{equation}
Note that the IMQ-based schemes are algebraically simple and closely parallel to the MQ and Gaussian constructions. Again, the IMQ forward scheme \eqref{eq:IMQ_Forward} is explicit, the IMQ backward scheme \eqref{eq:IMQ_Backward} is implicit, and the IMQ CN form \eqref{eq:IMQ_CN} provides a symmetric combination; in each case, the classical scheme can be recovered in the limit $\epsilon\to 0$.

We summarize the developed RBF-based CN scheme as
\begin{equation}\label{eq:RBF_CN}
    u_{n+1} = u_n + a_n \left( \mathcal{L}u_n + \mathcal{L}u_{n+1} \right),
\end{equation}
where the coefficient $a_n$ is given by
\begin{equation}
a_n =
\begin{cases}
    \displaystyle \frac{e^{\epsilon^2_n \Delta t^2} - 1}{2 \epsilon^2_n \Delta t} & \text{for Gaussian RBF}, \\[1.0em]
    \displaystyle \frac{\Delta t \,\sqrt{1 + \epsilon_n^2 \Delta t^2}}{1 + \sqrt{1 + \epsilon_n^2 \Delta t^2}} \, & \text{for MQ RBF}, \\[1.0em]
    \displaystyle \frac{\Delta t}{1 + \sqrt{1 + \epsilon_n^2 \Delta t^2}} \, & \text{for IMQ RBF}.
\end{cases}
\end{equation}
which can be written the compact form as 
\begin{equation}\label{eq:Compact_CN}
    \left( 1 - a_n \mathcal{L}\right)\, u_{n+1} = \left( 1 + a_n \mathcal{L}\right)\, u_{n}
\end{equation}

\section{Consistency}\label{sec:consistency}
We perform the consistency analysis focusing solely on the temporal aspect of the proposed RBF-CN schemes. The spatial operator $\mathcal{L}$ is regarded as an exact continuous linear operator without introducing any spatial discretization. Let $u(\mathbf{x},t)$ denote the exact solution of the continuous problem, and define
\begin{equation}
u_n := u(\mathbf{x}, t_n), \quad u_{n+1} := u(\mathbf{x}, t_{n+1}),
\end{equation}
as the exact solution values at discrete time levels $t_n$ and $t_{n+1}$, respectively. The local truncation error at time level $t_{n+1}$ is obtained by substituting the exact solution into the numerical scheme

\begin{equation}\label{eq:trunc_err1}
\tau_n = \frac{(1 - a_n \mathcal{L})\,u_{n+1} - (1 + a_n \mathcal{L})\,u_n}{\Delta t}.
\end{equation}

Expanding $u_{n+1}$ via Taylor series about $t_n$ and substituting $\frac{\partial^k u}{\partial t^k} = \mathcal{L}^k u$ by the linearity of $\mathcal{L}$, yields
\begin{equation}\label{eq:exp1}
u_{n+1} = \sum_{k=0}^{6} \frac{\Delta t^k}{k!} \, \mathcal{L}^k u_n + \mathcal{O}(\Delta t^7).
\end{equation}
Substituting this into the definition of $\tau_n$ yields
\begin{align}\label{eq:trunc_err2}
    \tau_n \,\Delta t & =   \left( 1 - a_n \mathcal{L}\right)\, \left(  \sum_{k=0}^6 \frac{\Delta t^k}{k!} \mathcal{L}^k u_n + \mathcal{O}(\Delta t^7) \right) - \left( 1 + a_n \mathcal{L}\right)\, u_n \nonumber \\
    & = - \,2 a_n \,\mathcal{L}u_n +  \left( 1 - a_n \mathcal{L}\right)\,   \sum_{k=1}^6 \frac{\Delta t^k}{k!} \mathcal{L}^k u_n + \mathcal{O}(\Delta t^7)  \nonumber \\
    & = \sum_{k=1}^6 \frac{\Delta t^k}{k!} \mathcal{L}^k u_n - a_n \left( \mathcal{L} u_n + \sum_{k=0}^6 \frac{\Delta t^k}{k!} \mathcal{L}^{k+1} u_n\right)  + \mathcal{O}(\Delta t^7).
\end{align}

\subsection{Gaussian RBF}
For the Gaussian RBF case, the coefficient $a_n$ admits the expansion
\begin{equation}\label{eq:Gaussian_a}
    a_n = \frac{\Delta t}{2} + \frac{\epsilon_n^2 \Delta t^3}{4} + \frac{\epsilon_n^4 \Delta t^5}{12}   + \mathcal{O}(\Delta t^7).
\end{equation}
using this \eqref{eq:Gaussian_a} in the expression \eqref{eq:trunc_err2} for local truncation error $\tau_n$, we have
\begin{align}\label{eq:trunc_err_GS1}
    \tau_n
    & =  \left( \mathcal{L} u_n - \frac{1}{2} \mathcal{L} u_n - \frac{1}{2}\mathcal{L} u_n \right)  + \Delta t \left( \frac{1}{2} \mathcal{L}^2u_n - \frac{1}{2} \mathcal{L}^2u_n\right) \nonumber \\
    & \quad + \Delta t^2 \left( \frac{1}{6} \mathcal{L}^3u_n - \frac{1}{4} \epsilon_n^2\mathcal{L}u_n - \frac{1}{4} \epsilon_n^2\mathcal{L}u_n - \frac{1}{4} \mathcal{L}^3u_n\right) \nonumber\\
    & \quad + \Delta t^3 \left( \frac{1}{24} \mathcal{L}^4u_n - \frac{1}{12} \mathcal{L}^4u_n - \frac{1}{4} \epsilon_n^2\mathcal{L}^2u_n \right) \\
    & \quad + \Delta t^4 \left( \frac{1}{120} \mathcal{L}^5u_n - \frac{1}{48} \mathcal{L}^5u_n - \frac{1}{8} \epsilon_n^2\mathcal{L}^3u_n - \frac{1}{12} \epsilon_n^4 \mathcal{L}u_n - \frac{1}{12} \epsilon_n^4 \mathcal{L}u_n\right) \nonumber\\
    & \quad + \Delta t^5 \left( \frac{1}{720} \mathcal{L}^6u_n - \frac{1}{240} \mathcal{L}^6u_n - \frac{1}{24} \epsilon_n^2\mathcal{L}^4u_n - \frac{1}{12} \epsilon_n^4 \mathcal{L}^2u_n \right) + \mathcal{O}(\Delta t^6) \nonumber
\end{align}
Further simplification leads to 
\begin{align}\label{eq:trunc_err_GS2}
    \tau_n 
    & =  - \frac{\Delta t^2}{12} \mathcal{L}\left( \mathcal{L}^2u_n +  6 \epsilon_n^2u_n \right)  - \frac{\Delta t^3}{24} \mathcal{L}^2\left( \mathcal{L}^2u_n +  6 \epsilon_n^2 u_n \right) \nonumber\\
    & \quad - \frac{\Delta t^4}{240} \mathcal{L}\left( 3\mathcal{L}^4u_n +  30 \epsilon_n^2\mathcal{L}^2 u_n + 40  \epsilon_n^4 u_n\right) \\
    & \quad - \frac{\Delta t^5}{360} \mathcal{L}^2 \left( \mathcal{L}^4u_n +  15 \epsilon_n^2\mathcal{L}^2 u_n + 30 \epsilon_n^4 u_n\right) + \mathcal{O}(\Delta t^6) \nonumber
\end{align}

From \eqref{eq:trunc_err_GS2}, it is evident that the Gaussian RBF-CN scheme retains the $\mathcal{O}(\Delta t^2)$ local truncation error characteristic of the standard CN method, thereby preserving its inherent second–order temporal accuracy. However, the coefficients of the higher–order terms ($\Delta t^2$ onward) are not uniquely determined, as they depend on the choice of the shape parameter $\epsilon_n$. By appropriately tuning $\epsilon_n^2$, these coefficients can be successively nullified, thereby elevating the scheme to higher–order temporal accuracy.

\subsubsection*{Third/Fourth Order Accuracy:}  
The scheme attains third–order temporal accuracy if the shape parameter $\epsilon_n^2$ is chosen such that the $\Delta t^2$ term in the truncation error vanishes, namely  
\begin{equation}\label{eq:GS_ep34_1}
    \mathcal{L}^2u_n + 6\,\epsilon_n^2\,u_n = 0.
\end{equation}
Expressing in terms of temporal derivatives yields
\begin{equation}\label{eq:GS_eps34_2}
    \epsilon_n^2 = -\frac{1}{6u_n}\,\frac{\partial^2 u_n}{\partial t^2},
\end{equation}
Interestingly, the expression for $\tau_n$ in \eqref{eq:trunc_err_GS2} reveals that the $\Delta t^2$ and $\Delta t^3$ terms have the same operator structure. Thereby, the same choice of $\epsilon_n^2$ also cancels the $\Delta t^3$ term, elevating the scheme to fourth–order accuracy without any additional constraints.

\subsubsection*{Fifth Order Accuracy:} The fifth order accuracy can be achieved if all the terms up to order of $\Delta t^4$ vanish in the truncation error expression, which implies
\begin{equation}\label{eq:GS5_1}
      20 \left( \mathcal{L}^2u_n +  6 \epsilon_n^2u_n \right)  + 10 \, \Delta t \left( \mathcal{L}^3u_n +  6 \epsilon_n^2 \mathcal{L}u_n \right) + \Delta t^2 \left( 3\mathcal{L}^4 u_n +  30 \epsilon_n^2\mathcal{L}^2 u_n + 40  \epsilon_n^4 u_n\right)  = 0 .
\end{equation}
The two possible values of $\epsilon^2_n$ from the expression \eqref{eq:GS5_1} are 
\begin{equation}\label{eq:GS_eps5}
    \epsilon_{n}^{2,\pm} 
    = \frac{-B \,\pm\, \sqrt{B^{2} - 4AC}}{2A},
\end{equation}
where
\begin{align*}
    A &= 40\,\Delta t^{2} \, u_{n}, \\
    B &= 30\left( \Delta t^{2} \mathcal{L}^{2} u_{n} \,+\, 2\Delta t \,\mathcal{L} u_{n} \,+\, 4u_{n} \right), \\
    C &= 3\Delta t^{2} \mathcal{L}^{4} u_{n} \,+\, 10\Delta t \,\mathcal{L}^{3} u_{n} \,+\, 20 \mathcal{L}^{2} u_{n}.
\end{align*}
We next demonstrate that at least one of the $\epsilon_{n}^{2, \pm}$ values yields consistency of the scheme. Assume that $u_n \neq 0$. First, consider the case $u_n > 0$:
\begin{align*}
    \lim_{\Delta t \to 0} \epsilon_{n}^{2,+} 
    & = \lim_{\Delta t \to 0}\frac{-B + \sqrt{B^2 - 4AC} }{2A}  \\
    & = -\lim_{\Delta t \to 0}\frac{2C}{ B + \sqrt{B^2 - 4AC}} \\
    & = -2 \frac{ \displaystyle\lim_{\Delta t \to 0} C}{ \displaystyle\lim_{\Delta t \to 0} \left( B + \sqrt{B^2 - 4AC} \right)} \\
    & = -\frac{\mathcal{L}^2 u_n}{3 \left( u_n + |u_n| \right)} 
      = -\frac{1}{6u_n}\,\frac{\partial^2 u_n}{\partial t^2},
\end{align*}
which corresponds to the optimal value of $\epsilon_n^{2}$ for achieving fourth-order accuracy of the scheme. This implies that $\tau_n \to 0$ as $\Delta t \to 0$, and therefore consistency holds for $\epsilon_n^{2,+}$. Similarly, consider the limit
\begin{align*}
   \lim_{\Delta t \to 0} \Delta t^2 \,\epsilon_{n}^{2,-}  
   & = \lim_{\Delta t \to 0}\frac{ \Delta t^2 \left(-B - \sqrt{B^2 - 4AC} \right) }{2A}  \\
   &  = \frac{  \displaystyle\lim_{\Delta t \to 0} \left(-B - \sqrt{B^2 - 4AC} \right) }{80 u_n}  \\
   & = \frac{ 120  \left(- u_n - |u_n|  \right) }{80 u_n} = -3,
\end{align*}
Thus, $|\epsilon_n^{2,-}| \le \frac{C}{\Delta t^2}$, for some constant $C>0$. To examine the consistency, consider the third term of \eqref{eq:trunc_err_GS2}. We have
\begin{align*}
    \left | \frac{\Delta t^4}{240} \mathcal{L}\left( 3\mathcal{L}^4u_n +  30 \epsilon_n^2\mathcal{L}^2 u_n + 40  \epsilon_n^4 u_n\right)\right|
    & \le \frac{\Delta t^4}{80} |\mathcal{L}^5u_n| + \frac{\Delta t^4}{8} |\mathcal{L}^3u_n| |\epsilon_n^2| + \frac{\Delta t^4}{6} |\mathcal{L}u_n| |\epsilon_n^4| \\
    & \le \frac{\Delta t^4}{80} |\mathcal{L}^5u_n| + \frac{\Delta t^2}{8} |\mathcal{L}^3u_n| C + \frac{1}{6} |\mathcal{L}u_n| C^2.
\end{align*}
Clearly, the term does not vanish as $\Delta t \to 0$. Thus, the consistency does not hold for $\epsilon_n^{2,-}$ when $u_n > 0$. By repeating the similar analysis for the case when $u_n <0$, the consistency only hold for $\epsilon_n^{2,-}$ and fails when we choose $\epsilon_n^{2,+}$. This indicates that, among the two roots $\epsilon_n^{2,\pm}$, only one ensures consistency, which depends on the sign of $u_n$. The other root grows unbounded as $\Delta t \to 0$ and does not yield a consistent scheme.

In a similar way, one can find the optimal shape parameter $\epsilon_n^2$ to achieve even higher accuracy, which will not be discussed here.

\begin{remark}
It is worth noting that the denominator in all forms of the optimal $\epsilon_n^2$ expressions (see \eqref{eq:GS_eps34_2} and \eqref{eq:GS_eps5}) vanishes when $u_n = 0$, causing $|\epsilon_n^2|$ to approach infinity. In such cases, a sufficiently large finite value can be assigned to $\epsilon_n^2$ to avoid numerical instability, although this may slightly compromise the optimal order of accuracy.
\end{remark}

\subsection{MQ RBF}
For the MQ RBF-based CN scheme, the coefficient $a_n$ is given by
\begin{equation}\label{eq:MQ_a_series}
    a_n = \frac{\Delta t}{2} + \frac{\epsilon_n^2 \Delta t^3}{8} - \frac{\epsilon_n^4 \Delta t^5}{16}  + \mathcal{O}(\Delta t^7).
\end{equation}
Substituting \eqref{eq:MQ_a_series} into the general truncation error expression yields
\begin{align}\label{eq:trunc_err_MQ1}
    \tau_n
    & =  \left( \mathcal{L} u_n - \frac{1}{2} \mathcal{L} u_n - \frac{1}{2}\mathcal{L} u_n \right)  
      + \Delta t \left( \frac{1}{2} \mathcal{L}^2u_n - \frac{1}{2} \mathcal{L}^2u_n\right) \nonumber \\
    & \quad + \Delta t^2 \left( \frac{1}{6} \mathcal{L}^3u_n - \frac{1}{8} \epsilon_n^2\mathcal{L}u_n - \frac{1}{8} \epsilon_n^2\mathcal{L}u_n - \frac{1}{4} \mathcal{L}^3u_n\right) \nonumber\\
    & \quad + \Delta t^3 \left( \frac{1}{24} \mathcal{L}^4u_n - \frac{1}{12} \mathcal{L}^4u_n - \frac{1}{8} \epsilon_n^2\mathcal{L}^2u_n \right) \\
    & \quad + \Delta t^4 \left( \frac{1}{120} \mathcal{L}^5u_n - \frac{1}{48} \mathcal{L}^5u_n - \frac{1}{16} \epsilon_n^2\mathcal{L}^3u_n + \frac{1}{16} \epsilon_n^4 \mathcal{L}u_n + \frac{1}{16} \epsilon_n^4 \mathcal{L}u_n\right) \nonumber\\
    & \quad + \Delta t^5 \left( \frac{1}{720} \mathcal{L}^6u_n - \frac{1}{240} \mathcal{L}^6u_n - \frac{1}{48} \epsilon_n^2\mathcal{L}^4u_n + \frac{1}{16} \epsilon_n^4 \mathcal{L}^2u_n \right) 
      + \mathcal{O}(\Delta t^6) \nonumber
\end{align}

This further simplifies to
\begin{align}\label{eq:trunc_err_MQ2}
    \tau_n
    & = -\frac{\Delta t^2}{12}\,\mathcal{L}\left(\mathcal{L}^2 u_n + 3\epsilon_n^2 u_n\right)
      -\frac{\Delta t^3}{24}\,\mathcal{L}^2\left(\mathcal{L}^2 u_n + 3\epsilon_n^2 u_n\right) \nonumber\\
    &\quad -\frac{\Delta t^4}{80}\,\mathcal{L}\left(\mathcal{L}^4 u_n + 5\epsilon_n^2\mathcal{L}^2 u_n - 10\epsilon_n^4 u_n\right) \\
    &\quad -\frac{\Delta t^5}{720}\,\mathcal{L}^2\left(2\mathcal{L}^4 u_n + 15\epsilon_n^2\mathcal{L}^2 u_n - 45\epsilon_n^4 u_n\right)
      + \mathcal{O}(\Delta t^6) \nonumber
\end{align}

\subsubsection*{Third/Fourth Order Accuracy:}  
For the MQ–RBF case, elimination of the $\Delta t^2$ term in \eqref{eq:trunc_err_MQ2} requires  
\begin{equation}\label{eq:MQ_eps34_1}
    \mathcal{L}^2 u_n + 3\,\epsilon_n^2\,u_n = 0,
\end{equation}
which is equivalent to  
\begin{equation}\label{eq:MQ_eps34_2}
    \epsilon_n^2 = -\frac{1}{3u_n}\,\frac{\partial^2 u_n}{\partial t^2}.
\end{equation}
As in the Gaussian case, this same choice of $\epsilon_n^2$ also removes the $\Delta t^3$ term, yielding fourth–order accuracy without further conditions.

\subsubsection*{Fifth Order Accuracy:}
For the MQ-based scheme, fifth-order temporal accuracy is obtained when
\begin{equation}\label{eq:MQ5_1}
      20\left(\mathcal{L}^2 u_n + 3\epsilon_n^2 u_n\right) + 10 \Delta t\left(\mathcal{L}^3 u_n + 3\epsilon_n^2 \mathcal{L}u_n\right) + 3 \Delta t^2\left(\mathcal{L}^4 u_n + 5\epsilon_n^2\mathcal{L}^2 u_n - 10\epsilon_n^4 u_n\right)  = 0 .
\end{equation}
This condition yields a quadratic equation for $\epsilon_n^2$ with roots
\begin{equation}\label{eq:MQ_eps5}
    \epsilon_{n}^{2,\pm} 
    = \frac{-B \,\pm\, \sqrt{B^{2} - 4AC}}{2A},
\end{equation}
where
\begin{align*}
    A &= -30\,\Delta t^{2} \, u_{n}, \\
    B &= 15\left( \Delta t^{2} \mathcal{L}^{2} u_{n} \,+\, 2\Delta t \,\mathcal{L} u_{n} \,+\, 4u_{n} \right), \\
    C &= 3\Delta t^{2} \mathcal{L}^{4} u_{n} \,+\, 10\Delta t \,\mathcal{L}^{3} u_{n} \,+\, 20 \mathcal{L}^{2} u_{n}.
\end{align*}
As like Gaussian case, consider the case $u_n > 0$. Taking $\Delta t \to 0$ in \eqref{eq:MQ_eps5} gives
\begin{align*}
    \lim_{\Delta t \to 0} \epsilon_{n}^{2,+} 
    & = -\frac{2 \mathcal{L}^2 u_n}{3 \left( u_n + |u_n| \right)} = -\frac{1}{3u_n}\,\frac{\partial^2 u_n}{\partial t^2},
\end{align*}
This matches the optimal $\epsilon_n^2$ for fourth-order accuracy, so $\epsilon_{n}^{2,+}$ ensures consistency. For the other branch,
\begin{align*}
   \lim_{\Delta t \to 0} \Delta t^2 \,\epsilon_{n}^{2,-}  
   & = \lim_{\Delta t \to 0}\frac{ \Delta t^2 \left(-B - \sqrt{B^2 - 4AC} \right) }{2A}  = \frac{  \left( u_n + |u_n|  \right) }{ u_n} = 2,
\end{align*}
which implies $\epsilon_n^{2,-} = \mathcal{O}(\Delta t^{-2})$.
implies $\epsilon_n^{2,-} = \mathcal{O}\left(\frac{1}{\Delta t^2} \right)$. Thus, by choosing $|\epsilon_n^{2,-}| \le \frac{C}{\Delta t^2}$, for some constant $C>0$, one can show that the third term of \eqref{eq:trunc_err_MQ2} remain finite as $\Delta t \to 0$, so the truncation error does not vanish. When $u_n < 0$, the argument is analogous, but the consistent choice is $\epsilon_{n}^{2,-}$. 

\subsection{IMQ RBF}
For the IMQ RBF-based CN scheme, the expansion of the coefficient $a_n$ is given by
\begin{equation}\label{eq:IMQ_a_series}
    a_n = \frac{\Delta t}{2} - \frac{\epsilon_n^2 \Delta t^3}{8} + \frac{\epsilon_n^4 \Delta t^5}{16} + \mathcal{O}(\Delta t^7).
\end{equation}
Inserting into the expression for $\tau_n$ and simplifying yields 
\begin{align}\label{eq:trunc_err_IMQ2}
    \tau_n
    &= -\frac{\Delta t^2}{12}\,\mathcal{L}\left(\mathcal{L}^2 u_n - 3\epsilon_n^2 u_n \right)
      -\frac{\Delta t^3}{24}\,\mathcal{L}^2\left(\mathcal{L}^2 u_n - 3\epsilon_n^2 u_n \right) \nonumber\\
    &\quad -\frac{\Delta t^4}{80}\,\mathcal{L}\left(\mathcal{L}^4 u_n - 5\epsilon_n^2\mathcal{L}^2 u_n + 10\epsilon_n^4 u_n \right) \\
    &\quad -\frac{\Delta t^5}{720}\,\mathcal{L}^2\left(2\mathcal{L}^4 u_n - 15\epsilon_n^2\mathcal{L}^2 u_n + 45\epsilon_n^4 u_n \right)
      + \mathcal{O}(\Delta t^6) \nonumber
\end{align}

\subsubsection*{Third/Fourth Order Accuracy:}  
For the IMQ-RBF scheme, the optimal shape parameter for third and fourth order temporal accuracy is  
\begin{equation}\label{eq:IMQ_eps34}
\epsilon_n^2 = \frac{1}{3u_n}\frac{\partial^2 u_n}{\partial t^2}.
\end{equation}

\subsubsection*{Fifth Order Accuracy:}
The optimal $\epsilon_n^2$ for achieving fifth-order temporal accuracy is
\begin{equation}\label{eq:IMQ_eps5}
    \epsilon_{n}^{2,\pm} 
    = \frac{-B \,\pm\, \sqrt{B^{2} - 4AC}}{2A},
\end{equation}
where
\begin{align*}
    A &= 30\,\Delta t^{2} \, u_{n}, \\
    B &= -15\left( \Delta t^{2} \mathcal{L}^{2} u_{n} \,+\, 2\Delta t \,\mathcal{L} u_{n} \,+\, 4u_{n} \right), \\
    C &= 3\Delta t^{2} \mathcal{L}^{4} u_{n} \,+\, 10\Delta t \,\mathcal{L}^{3} u_{n} \,+\, 20 \mathcal{L}^{2} u_{n}.
\end{align*}
It can be shown that the consistent branch of \eqref{eq:IMQ_eps5} depends on the sign of $u_n$. For $u_n > 0$, the appropriate choice is $\epsilon_{n}^{2,-}$, while for $u_n < 0$, it is $\epsilon_{n}^{2,+}$. This behavior is the reverse of that observed for the Gaussian and MQ cases.

\begin{remark}
In all RBF cases, the asymptotic expansion of $a_n$ is derived under the assumption $|\epsilon_n^2|\,\Delta t^2 < 1$. This guarantees convergence of the series, ensures that $a_n$ is real-valued and yields $a_n>0$.
\end{remark}

\section{Stability Analysis}\label{sec:stability}
We analyze the stability of the RBF-based CN scheme by considering the evolution of a single Fourier mode. Let the numerical solution be represented as a superposition of modes of the form
\begin{equation*}
u_j^n = \hat{u}^n e^{i k j h},
\end{equation*}
where $k$ is the wave number, $j$ indexes the spatial points, and $h$ is the spatial step size. Substituting this mode into \eqref{eq:Compact_CN} yields
\begin{equation*}
\left( 1 - a_n \lambda_k \right) \hat{u}^{\,n+1} 
= \left( 1 + a_n \lambda_k \right) \hat{u}^{\,n},
\end{equation*}
where $\lambda_k$ denotes the eigenvalue of the discrete spatial operator $\mathcal{L}$ corresponding to mode $k$.  
Defining the amplification factor $G_k$ as
\begin{equation*}
\hat{u}^{n+1} = G_k \, \hat{u}^n, \qquad 
G_k = \frac{1 + a_n \lambda_k}{1 - a_n \lambda_k}, 
\end{equation*} 
the stability condition requires $|G_k| \le 1$ for all eigenvalues $\lambda_k$.

\subsection{General stability condition}
Let $\lambda_k \in \mathbb{C}$ be expressed as $\lambda_k = \lambda_r + i \lambda_c$. The stability requirement $|G_k| \le 1$ is equivalent to
\begin{align}
\left|\frac{1 + a_n \lambda_k}{1 - a_n \lambda_k}\right| \le 1 
&\Longleftrightarrow (1 + a_n \lambda_r)^2 + (a_n \lambda_c)^2 \le (1 - a_n \lambda_r)^2 + (a_n \lambda_c)^2 \nonumber \\
&\Longleftrightarrow a_n \, \Re(\lambda_k) \le 0 \quad \forall k. \label{eq:stab_gen}
\end{align} 
For the RBF-based schemes considered here, $a_n$ is real and strictly positive by construction, provided the shape parameter satisfies
\begin{equation}\label{eq:eps_dt_condition}
|\epsilon_n^2|\,\Delta t^2 < 1.
\end{equation}
Under this assumption, the stability condition simplifies to
\begin{equation}\label{eq:stab}
\Re(\lambda_k) \le 0 \quad \forall k.
\end{equation}

\begin{remark}[Pure Diffusion]
If $\mathcal{L}$ is the standard diffusion operator discretized using second-order central differences on a uniform grid, then $\mathcal{L}$ is self-adjoint and negative semi-definite in the discrete $\ell^2$ inner product. All eigenvalues satisfy $\lambda_k \le 0$, and the above condition is automatically satisfied. In this case, the RBF-based CN scheme is unconditionally stable.
\end{remark}

\begin{remark}[Convection--Diffusion]
For a convection--diffusion operator discretized with central differences, $\mathcal{L}$ is generally non-self-adjoint and may have complex eigenvalues.  

When diffusion dominates, one typically has $\Re(\lambda_k) \le 0$ ensuring stability under the present scheme.  

In the pure convection case ($\Re(\lambda_k) = 0$), the scheme is neutrally stable ($|G_k| = 1$) and introduces no numerical damping.  

When convection and diffusion are of comparable strength, $\mathcal{L}$ may have modes with $\Re(\lambda_k) > 0$, violating the stability condition and leading to spectral instability. Even when $\Re(\lambda_k) \le 0$ for all $k$, strong non-normality can cause transient growth (pseudospectral effects), so eigenvalue checks alone may not guarantee robust stability.
\end{remark}

For convenience, we summarize in Tables \ref{tab:RBF_CN_summarya} and \ref{tab:RBF_CN_summaryb} the optimal shape parameters $\epsilon_n^2$ for all RBF schemes considered. These tables represent the results previously discussed in the individual subsections and indicates the temporal order of accuracy for which each $\epsilon_n^2$ is optimal. 

\begin{sidewaystable*}
\caption{\label{tab:RBF_CN_summarya}Summary of Standard and RBF-based CN schemes.}
\begin{ruledtabular}
\begin{tabular}{cccc}
& & Order of accuracy & Optimal shape parameter \\ \hline
Standard CN & $\displaystyle u_{n+1} = u_n + \frac{\Delta t}{2} (\mathcal{L} u_n + \mathcal{L} u_{n+1})$ & $\mathcal{O}(\Delta t^2)$ & --- \\
Gaussian RBF-CN & $u_{n+1} = u_n + a_n (\mathcal{L} u_n + \mathcal{L} u_{n+1})$ & $\mathcal{O}(\Delta t^4)$ & $\displaystyle -\frac{1}{6 u_n} \frac{\partial^2 u_n}{\partial t^2}$ \\
& $\displaystyle a_n = \frac{e^{\epsilon_n^2 \Delta t^2}-1}{2 \epsilon_n^2 \Delta t}$ & & \\ 
& & O($\Delta t^5$) & 
$\epsilon_n^2 = 
\begin{cases} 
\epsilon_n^{2,+}, & u_n>0 \\ 
\epsilon_n^{2,-}, & u_n<0 
\end{cases}$ \\
& & & $\displaystyle \epsilon_n^{2,\pm} = \frac{-B \pm \sqrt{B^2 - 4AC}}{2A}$ \\
& & & $A = 40 \Delta t^2 u_n$ \\
& & & $B = 30 (\Delta t^2 \mathcal{L}^2 u_n + 2\Delta t \mathcal{L} u_n + 4 u_n)$ \\
& & & $\displaystyle C = 3 \Delta t^2 \mathcal{L}^4 u_n + 10 \Delta t \mathcal{L}^3 u_n + 20 \mathcal{L}^2 u_n$ \\ 
MQ RBF-CN & $u_{n+1} = u_n + a_n (\mathcal{L} u_n + \mathcal{L} u_{n+1})$ & O($\Delta t^4$) & $\displaystyle -\frac{1}{3 u_n} \frac{\partial^2 u_n}{\partial t^2}$ \\
& $\displaystyle a_n = \frac{\Delta t \sqrt{1+\epsilon_n^2 \Delta t^2}}{1 + \sqrt{1+\epsilon_n^2 \Delta t^2}}$ & & \\
& & O($\Delta t^5$) & 
$\epsilon_n^2 = 
\begin{cases} 
\epsilon_n^{2,+}, & u_n>0 \\ 
\epsilon_n^{2,-}, & u_n<0 
\end{cases}$ \\
& & & $\displaystyle \epsilon_n^{2,\pm} = \frac{-B \pm \sqrt{B^2 - 4AC}}{2A}$ \\
& & & $A = -30 \Delta t^2 u_n$ \\
& & & $B = 15 (\Delta t^2 \mathcal{L}^2 u_n + 2\Delta t \mathcal{L} u_n + 4 u_n)$ \\
& & & $\displaystyle C = 3 \Delta t^2 \mathcal{L}^4 u_n + 10 \Delta t \mathcal{L}^3 u_n + 20 \mathcal{L}^2 u_n$ \\
\end{tabular}
\end{ruledtabular}
\end{sidewaystable*}

\begin{sidewaystable*}
\caption{\label{tab:RBF_CN_summaryb}Summary of Standard and RBF-based CN schemes.}
\begin{ruledtabular}
\begin{tabular}{cccc}
& & Order of accuracy & Optimal shape parameter \\ \hline
IMQ RBF-CN & $u_{n+1} = u_n + a_n (\mathcal{L} u_n + \mathcal{L} u_{n+1})$ & O($\Delta t^4$) & $\displaystyle \frac{1}{3 u_n} \frac{\partial^2 u_n}{\partial t^2}$ \\
& $\displaystyle a_n = \frac{\Delta t}{1 + \sqrt{1+\epsilon_n^2 \Delta t^2}}$ & & \\
& & O($\Delta t^5$) & 
$\epsilon_n^2 = 
\begin{cases} 
\epsilon_n^{2,-}, & u_n>0 \\ 
\epsilon_n^{2,+}, & u_n<0 
\end{cases}$ \\
& & & $\displaystyle \epsilon_n^{2,\pm} = \frac{-B \pm \sqrt{B^2 - 4AC}}{2A}$ \\
& & & $A = 30 \Delta t^2 u_n$ \\
& & & $B = -15 (\Delta t^2 \mathcal{L}^2 u_n + 2\Delta t \mathcal{L} u_n + 4 u_n)$ \\
& & & $C = 3 \Delta t^2 \mathcal{L}^4 u_n + 10 \Delta t \mathcal{L}^3 u_n + 20 \mathcal{L}^2 u_n$ \\
\end{tabular}
\end{ruledtabular}
\end{sidewaystable*}

\section{Numerical Experiments}\label{sec:numerical}
In this section, we investigate the accuracy and robustness of the proposed RBF-based CN schemes through several test problems. For simplicity, all test problems are posed in one spatial dimension. For all cases, Spatial derivatives are approximated using standard second-order central difference schemes, while the temporal discretization is handled by the proposed RBF-based CN method. The truncation error of the schemes can then be expressed as
\begin{equation}
    \tau_n = \mathcal{O}(\Delta t^p) + \mathcal{O}(\Delta x^2),
\end{equation}
where $p$ denotes the temporal order of the RBF-CN scheme. Since our focus is on verifying the temporal accuracy of the proposed RBF schemes, we choose the spatial mesh size according to $\Delta x = \Delta t^{p/2}$, so that the spatial discretization error does not dominate the overall truncation error. With this choice, the leading-order error reduces to $ \tau_n = \mathcal{O}(\Delta t^p)$.

The numerical simulation using the RBF-based CN scheme requires the shape parameter to be determined in advance for computing the solution at the subsequent time level. The expression for the optimal shape parameter involves derivatives of the solution variable with respect to time beyond the first order. Approximating these derivatives through backward differences necessitates the availability of solutions from previous time levels, which introduces a start-up difficulty. One possible approach is to utilize the exact solution to advance the computation until the stage at which $\ epsilon^2$ can be evaluated, after which the RBF-CN scheme can be applied. In practice, however, this option is rarely feasible since exact solutions are generally unavailable for most problems.

Another approach is, in the absence of an exact solution, one may instead employ the standard CN scheme during the initial steps and subsequently switch to the RBF-CN scheme. These pre-computation steps should be performed with a finer step size $\delta t < \Delta t$ to preserve the formal order of accuracy of the RBF–CN method. The choice of $\delta t$ must ensure that the local truncation error introduced by the CN scheme is of higher order than the global error of the RBF–CN scheme with $\Delta t$; otherwise, the overall convergence rate will be governed by the standard CN used for the startup stage. Another alternative is to replace the refined CN steps by a higher-order scheme, for example, by applying Richardson extrapolation to CN or by using another scheme of the same order $p$ as the RBF-CN.

\begin{remark}
    If we choose $\delta t = \Delta t / n_{0}$, where $n_{0}$ denotes the number of subdivisions of the working step $\Delta t$, then the required condition on $n_{0}$ follows from 
\begin{equation*}
E_{\text{CN}} \;\leq\; E_{\text{RBF-CN}} 
\;\;\;\Rightarrow\;\;\;
\left(\frac{\Delta t}{n_{0}}\right)^{2} \leq C\, \Delta t^{p},
\end{equation*}
and hence 
\begin{equation*}
n_{0} \;\geq\; C^{1/2}\, \Delta t^{\tfrac{2-p}{2}},
\end{equation*}
for some constant $C$ independent of $\Delta t$. In practice, choosing such a large $n_{0}$ at the initialization stage may entail significant computational cost.
\end{remark}

\subsection{Test Problem 1}
Consider the one-dimensional diffusion equation
\begin{equation}\label{eq:problem1}
    \frac{\partial u}{\partial t} = \frac{\partial^2 u}{\partial x^2}, 
    \qquad 0 \le x \le 1, \quad 0 < t \le 1,
\end{equation}
subject to the initial and boundary conditions
\begin{equation*}
    u(x,0) = \sin(\pi x), \qquad u(0,t) = u(1,t) = 0.
\end{equation*}
The exact solution of this problem is given by
\begin{equation}\label{eq:problem1_sol}
    u(x,t) = e^{-\pi^2 t} \sin(\pi x).
\end{equation}
We first choose the optimal $\epsilon^2$, which theoretically yields fourth-order temporal accuracy in the RBF-CN schemes. The expression for this optimal parameter involves the second-order time derivative of the solution $u$. Approximating this derivative with an $s$-th order backward difference requires the availability of the solution at the previous $s+2$ time levels beforehand. 

Consider the case when these previous steps' solutions are obtained from the exact solution \eqref{eq:problem1_sol}. Table \ref{tab:problem1_tab1} presents the global error and order of convergence of the RBF-CN schemes for six successively refined time steps $\Delta t$. The global error is defined as the discrete infinity norm over space and time of the difference between the exact and numerical solutions. The order of convergence at each step is estimated via a least-squares fit of the global errors from the current and all previous refined steps, which we can also call the global order.

\begin{table*}
\caption{\label{tab:problem1_tab1}Comparison of global error and convergence order of the RBF-CN schemes for various approximations of $\epsilon^2$ when starting stage solutions are obtained from the exact solution.}
\begin{ruledtabular}
\begin{tabular}{cccccccc}
& &\multicolumn{2}{c}{$\epsilon^2 \sim \mathcal{O}(\Delta t^1)$}&\multicolumn{2}{c}{$\epsilon^2 \sim \mathcal{O}(\Delta t^2)$}&\multicolumn{2}{c}{$\epsilon^2 \sim \mathcal{O}(\Delta t^3)$}\\
Scheme & $N_t$ & Global Error & Order & Global Error & Order & Global Error & Order \\ \hline
Gaussian RBF-CN   & 16  & 3.2607e-03 &  --   & 1.3973e-03 &  --   & 6.1869e-04 &  --   \\
 & 32  & 5.9758e-04 & 2.45  & 1.3821e-04 & 3.34  & 3.7438e-05 & 4.05  \\
 & 64  & 9.1029e-05 & 2.58  & 1.1366e-05 & 3.47  & 1.9644e-06 & 4.15  \\
 & 128 & 1.2652e-05 & 2.67  & 8.1923e-07 & 3.58  & 9.4092e-08 & 4.23  \\
 & 256 & 1.6665e-06 & 2.74  & 5.4991e-08 & 3.67  & 4.6021e-09 & 4.27  \\
 & 512 & 2.1435e-07 & 2.79  & 3.3430e-09 & 3.74  & 4.1724e-10 & 4.17  \\
\hline
 MQ RBF-CN  & 16  & 4.6071e-03 &   --   & 1.3866e-03 &  --   & 8.9648e-04 &  --   \\
 & 32  & 6.6289e-04 &  2.80  & 1.1895e-04 & 3.54  & 5.6970e-05 & 3.98  \\
 & 64  & 9.4920e-05 &  2.80  & 9.0534e-06 & 3.63  & 4.0699e-06 & 3.89  \\
 & 128 & 1.2894e-05 &  2.82  & 6.3004e-07 & 3.70  & 2.7148e-07 & 3.89  \\
 & 256 & 1.6815e-06 &  2.85  & 4.1538e-08 & 3.76  & 1.7414e-08 & 3.90  \\
 & 512 & 2.1358e-07 &  2.88  & 2.6438e-09 & 3.81  & 3.1673e-09 & 3.70  \\
\hline
IMQ RBF-CN & 16  & 2.7788e-03 &  --   & 1.4025e-03 &  --   & 5.1549e-04 &  --   \\
 & 32  & 5.6674e-04 & 2.29  & 1.4750e-04 & 3.25  & 2.8070e-05 & 4.20  \\
 & 64  & 8.9107e-05 & 2.48  & 1.2511e-05 & 3.40  & 9.2257e-07 & 4.56  \\
 & 128 & 1.2532e-05 & 2.60  & 9.1357e-07 & 3.53  & 5.6422e-09 & 5.44  \\
 & 256 & 1.6588e-06 & 2.69  & 6.1966e-08 & 3.63  & 1.8962e-09 & 4.84  \\
 & 512 & 2.0912e-07 & 2.76  & 1.2390e-08 & 3.47  & 6.3699e-09 & 3.73  \\
\end{tabular}
\end{ruledtabular}
\end{table*}

\begin{figure}[htbp!]
  \centering
    (a) \hspace{3.2 in} (b) \\ 
    \includegraphics[width=0.495\textwidth]{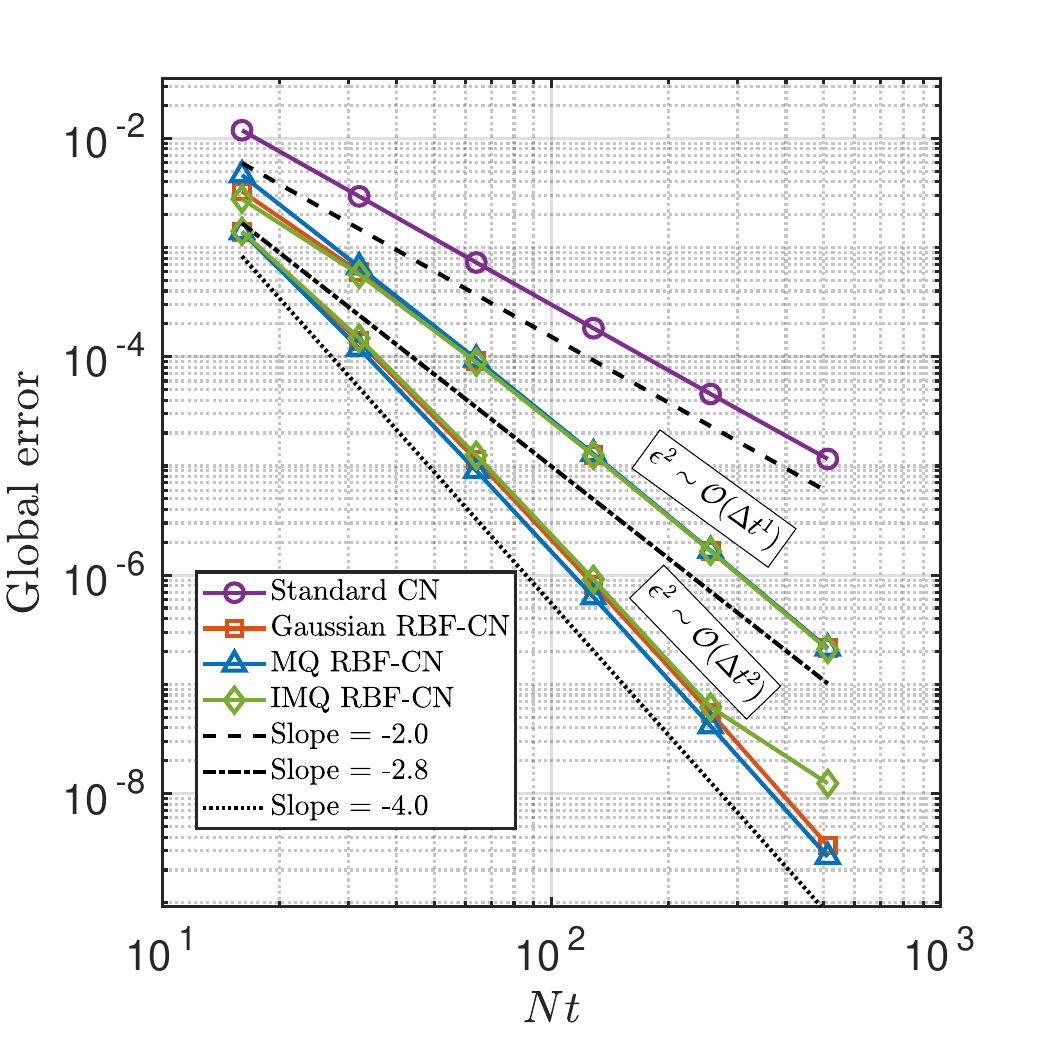}
    \includegraphics[width=0.495\textwidth]{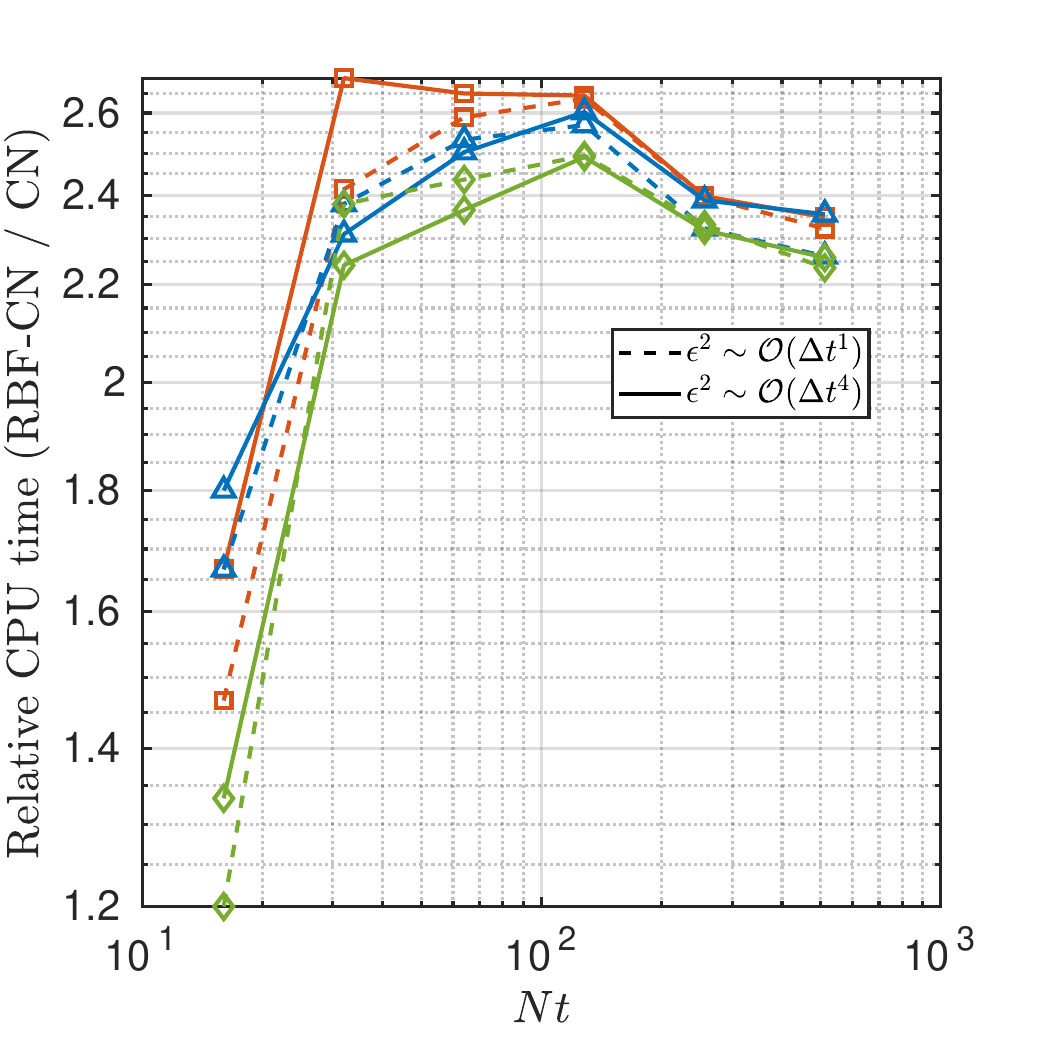}
  \caption{(a) Global convergence order and (b) CPU time of RBF-CN schemes relative to standard CN with optimal starting $\epsilon^2$ from the exact solution.}
  \label{fig:exactfig}
\end{figure}

Table \ref{tab:problem1_tab1} and Figure \ref{fig:exactfig}(a) demonstrate that the order of accuracy of the RBF-CN schemes is nearly $2.8$ for the first order approximation of $\epsilon^2$ with the finer $\Delta t$. The order improves close to the desired theoretical fourth order with the higher-order approximation of $\epsilon^2$. Additionally, observed that the MQ RBF-CN schemes provide relatively higher accuracy even for larger $\Delta t$, but the consistency-wise Gaussian CN shows better performance.

Although using the higher order approximation of $\epsilon^2$ in RBF-CN schemes improves the accuracy but it needs more starting steps to be computed beforehand. If these steps are taken from the exact solution, as shown in Figure \ref{fig:exactfig}(b), the CPU time consumed by RBF-CN schemes with any order of approximation of optimal $\epsilon^2$ is nearly the same, with a maximum of $2.6$ times that of CN. 

\begin{figure}[htbp!]
  \centering
    (a) \hspace{3.2 in} (b) \\ 
    \includegraphics[width=0.495\textwidth]{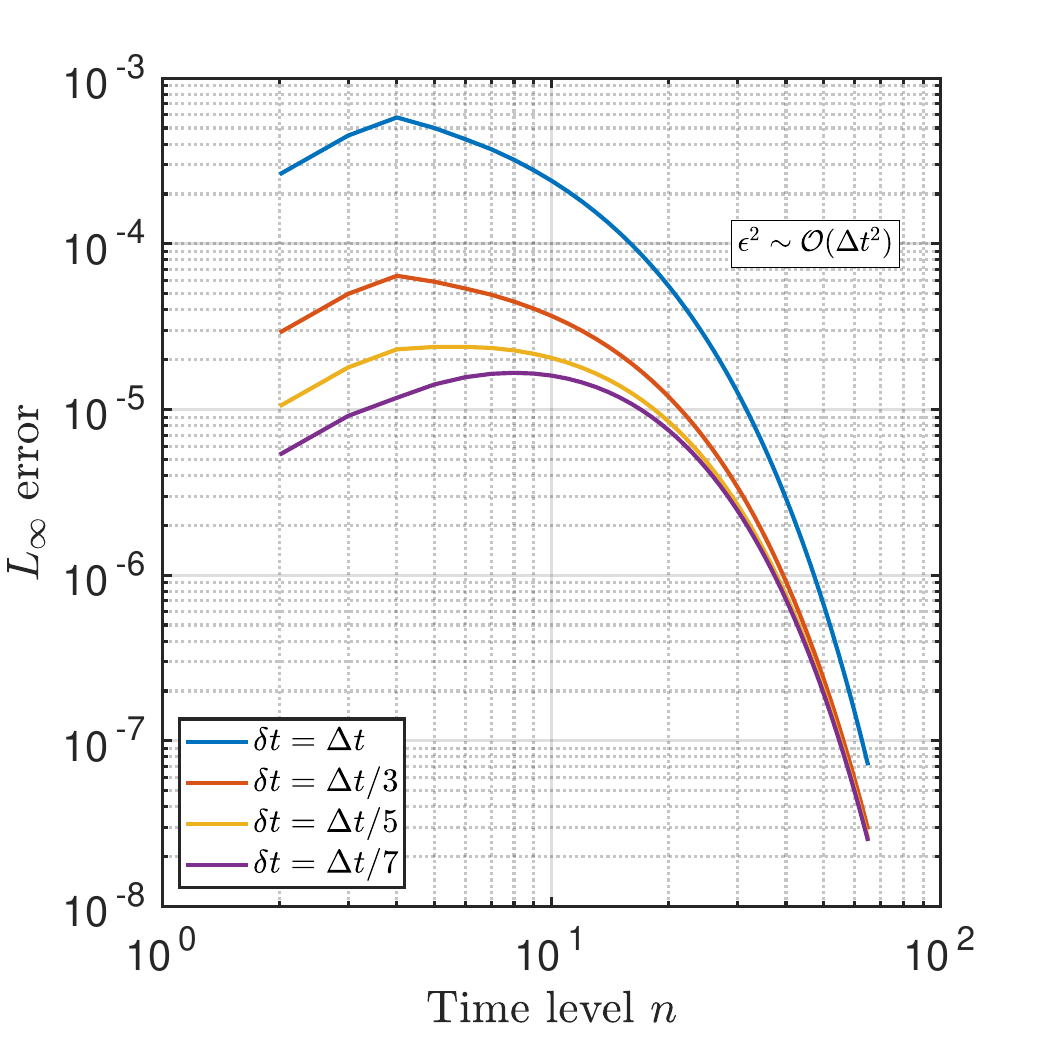}
    \includegraphics[width=0.495\textwidth]{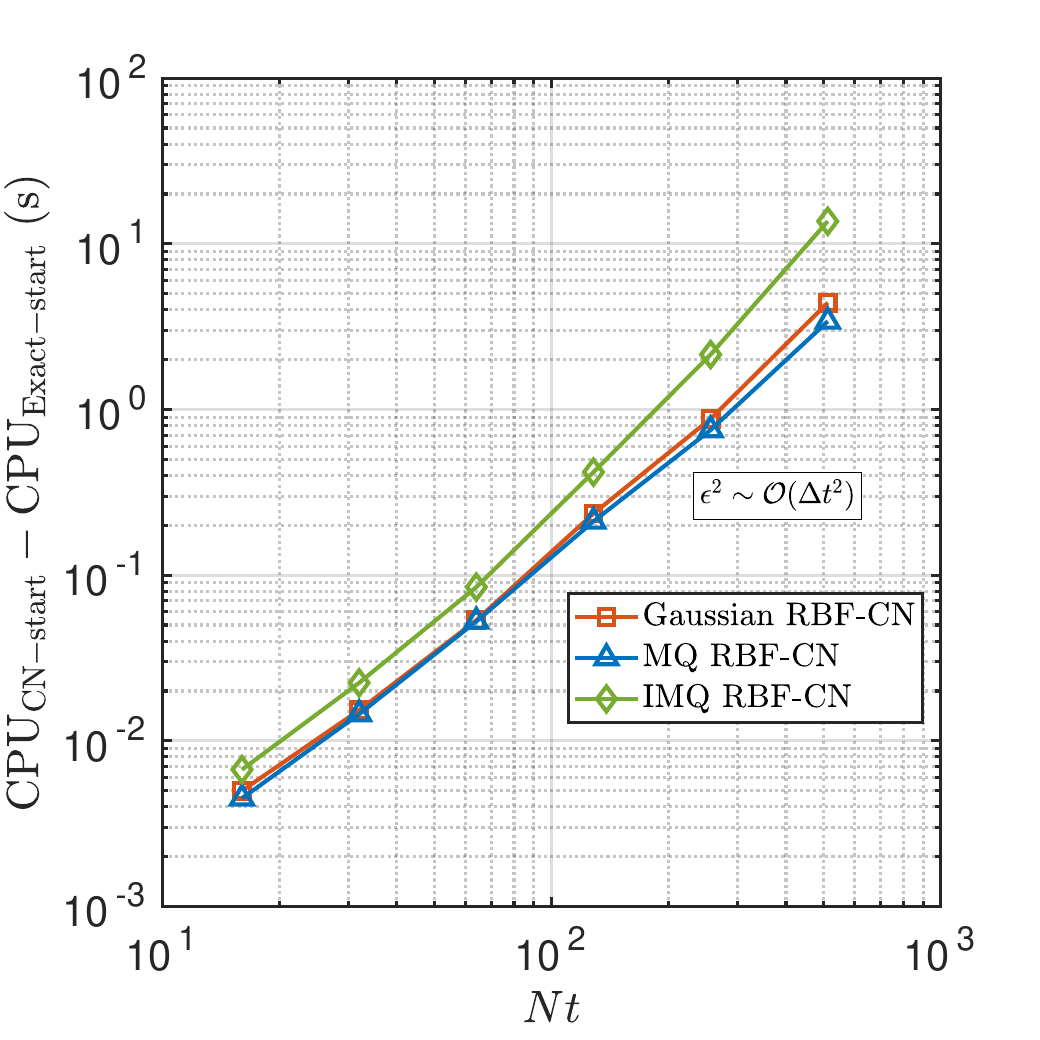}
  \caption{(a) Error using CN different $\delta t$ for the starting steps. (b) Additional CPU cost introduced by CN with $\delta t = \Delta t/40$ compared to the exact solution used for the starting steps.}
  \label{fig:CNstarting}
\end{figure}

We now examine the case where the initial steps are obtained numerically. When the standard CN scheme with $\delta t \gtrsim \Delta t/4$ is employed for this purpose, the global convergence order is observed to be governed by the CN scheme itself. This is because the leading error occurs in the starting steps, as demonstrated in Figure \ref{fig:CNstarting}(a). 

Table \ref{tab:problem1_tab2} presents the performance of the RBF-CN schemes when the starting steps are computed with the CN method using four different values of $\delta t$ smaller than $\Delta t/4$. The observed order is about three for $\delta t = \Delta t/5$ and improves with finer $\delta t$. At $\delta t = \Delta t/40$, the convergence order is nearly the same as when the exact solution is used to approximate the optimal $\epsilon^{2}$ (see the third column of Table \ref{tab:problem1_tab1}). For even finer $\delta t$, no further improvement in the order is observed. Achieving this accuracy requires roughly $120$ initial CN steps. Among the RBF choices, the MQ-based scheme performs better in terms of both accuracy and computational efficiency, whereas the IMQ-based scheme incurs a noticeably higher CPU time. To highlight this computational overhead, the additional CPU cost associated with employing the CN scheme for $\delta t = \Delta t/40$, relative to using the exact solution for the starting steps, is shown in Figure \ref{fig:CNstarting}(b). 

\begin{table*}
\caption{\label{tab:problem1_tab2} Comparison of convergence order and CPU time of the RBF-CN schemes with $\epsilon^2 \sim \mathcal{O}(\Delta t^2)$ when starting steps are computed using standard CN with various $\delta t$.}
\begin{ruledtabular}
\begin{tabular}{cccccccccc}
& &\multicolumn{2}{c}{RBF-CN Type I}&\multicolumn{2}{c}{RBF-CN Type I}&\multicolumn{2}{c}{RBF-CN Type I}&\multicolumn{2}{c}{RBF-CN Type I}\\
& & \multicolumn{2}{c}{$\delta t = \Delta t/5$}&\multicolumn{2}{c}{$\delta t = \Delta t/10$}&\multicolumn{2}{c}{$\delta t = \Delta t/20$}&\multicolumn{2}{c}{$\delta t = \Delta t/40$}\\
RBF & $N_t$ & Order & CPU Time & Order & CPU Time & Order & CPU Time & Order & CPU Time \\ \hline
Gaussian  & 16  &  --  & 0.0028  &  -- & 0.0036  &  -- &  0.0053 & --  & 0.0082 \\
 & 32  &  2.98 &  0.0131 & 3.25  & 0.0152  & 3.31  & 0.0190  &  3.33 & 0.0270  \\
 & 64  &  2.99 & 0.0893  & 3.35  & 0.0930  & 3.44  & 0.1074  &  3.46 & 0.1339  \\
 & 128 &  2.91 & 0.6909  &  3.40 & 0.7076  & 3.53  & 0.7665  & 3.57  & 0.8731  \\
 & 256 &  2.88 & 4.4969  & 3.38  & 4.5496  &  3.60 & 4.7793  & 3.65  &  5.1163 \\
 & 512 &  2.88 & 35.3404  & 3.36  & 35.5712  &  3.69 & 36.9044  & 3.74  &  37.9801 \\
\hline
 MQ & 16  &  --  & 0.0028  &  -- & 0.0036  &  -- &  0.0053 & --  & 0.0082 \\
 & 32  &  2.98 &  0.0131 & 3.25  & 0.0152  & 3.31  & 0.0190  &  3.33 & 0.0270  \\
 & 64  &  2.99 & 0.0893  & 3.35  & 0.0930  & 3.44  & 0.1074  &  3.46 & 0.1339  \\
 & 128 &  2.91 & 0.6909  &  3.40 & 0.7076  & 3.53  & 0.7665  & 3.57  & 0.8731  \\
 & 256 &  2.88 & 4.4969  & 3.38  & 4.5496  &  3.60 & 4.7793  & 3.65  &  5.1163 \\
 & 512 &  2.88 & 35.3404  & 3.36  & 35.5712  &  3.69 & 36.9044  & 3.74  &  37.9801 \\
\hline
 IMQ & 16  &  --  & 0.0038  &  --  & 0.0050  &  --  & 0.0063  &  --  & 0.0098 \\
 & 32  &  2.92 & 0.0169  & 3.17 & 0.0183  & 3.23 & 0.0238  & 3.24 & 0.0333 \\
 & 64  &  2.97 & 0.0989  & 3.29 & 0.1162  & 3.38 & 0.1380  & 3.40 & 0.1626 \\
 & 128 &  2.91 & 0.8098  & 3.37 & 0.8407  & 3.49 & 0.8662  & 3.52 & 1.0295 \\
 & 256 &  2.88 & 5.4020  & 3.37 & 5.5949  & 3.56 & 5.7734  & 3.61 & 6.2090 \\
 & 512 &  2.89 & 42.1665 & 3.37 & 41.7294 & 3.61 & 43.6573 & 3.67 & 46.0523 \\
\end{tabular}
\end{ruledtabular}
\end{table*}

Instead of employing the Crank–Nicolson (CN) scheme, a more efficient strategy is to use a method of the same fourth-order accuracy as the RBF-CN scheme for the initial steps. A straightforward choice is the standard CN scheme combined with Richardson extrapolation. Table \ref{tab:problem1_tab3} presents the global error, order of convergence, and CPU time for the RBF-CN scheme when the initial step computations are performed using CN with Richardson extrapolation. The results indicate that the convergence order closely follows the exact cases' values up to the time step size $\Delta = 1/256$, but decreases thereafter.
The decline may be because Richardson extrapolation relies on error cancellation, which becomes ineffective at very fine time steps where round-off and higher-order terms dominate. Nevertheless, with an appropriately chosen $\Delta t$, the CN scheme with Richardson extrapolation remains an effective choice for the initial stage. It requires only $9$ CN steps, leading to reduced computational cost, and the step count remains fixed for any problem, in contrast to the earlier approach of employing CN with progressively refined time steps. 

\begin{table*}
\caption{\label{tab:problem1_tab3} Comparison of global error, convergence order, and CPU time of the RBF-CN schemes with $\epsilon^2 \sim \mathcal{O}(\Delta t^2)$ when initial steps are computed using Richardson extrapolation and Gauss–Legendre IRK scheme.}
\begin{ruledtabular}
\begin{tabular}{cccccccc}
& &\multicolumn{3}{c}{RBF-CN Type II}&\multicolumn{3}{c}{RBF-CN Type III}\\
Scheme & $N_t$ & Global Error & Order & CPU Time
&Global Error & Order & CPU Time \\ \hline
Gaussian RBF-CN & 16  & 1.3631e-03 &  --  & 0.0028 & 1.3794e-03 &  --  & 0.0040 \\
 & 32  & 1.3428e-04 & 3.34 & 0.0126 & 1.3627e-04 & 3.34 & 0.0171 \\
 & 64  & 1.1190e-05 & 3.46 & 0.0891 & 1.1280e-05 & 3.47 & 0.1078 \\
 & 128 & 8.1182e-07 & 3.57 & 0.6833 & 8.1572e-07 & 3.58 & 0.8021 \\
 & 256 & 5.5856e-08 & 3.65 & 4.4915 & 5.4803e-08 & 3.66 & 4.9678 \\
 & 512 & 3.1573e-08 & 3.27 & 35.4310 & 3.2836e-09 & 3.74 & 36.6170 \\ \hline 
 MQ RBF-CN & 16  & 1.3523e-03 &  --  & 0.0027 & 1.3687e-03 &  --  & 0.0045 \\
 & 32  & 1.1502e-04 & 3.56 & 0.0123 & 1.1701e-04 & 3.55 & 0.0194 \\
 & 64  & 8.8775e-06 & 3.63 & 0.0850 & 8.9676e-06 & 3.63 & 0.1150 \\
 & 128 & 6.2265e-07 & 3.69 & 0.6433 & 6.2654e-07 & 3.70 & 0.7746 \\
 & 256 & 4.2465e-08 & 3.74 & 4.2401 & 4.1484e-08 & 3.76 & 4.7084 \\
 & 512 & 3.1573e-08 & 3.29 & 33.3956 & 2.6507e-09 & 3.80 & 34.6695 \\
\hline
 IMQ RBF-CN & 16  & 1.3682e-03 &  --  & 0.0028 & 1.3846e-03 &  --  & 0.0046 \\
 & 32  & 1.4536e-04 & 3.25 & 0.0122 & 1.4555e-04 & 3.25 & 0.0199 \\
 & 64  & 1.2334e-05 & 3.40 & 0.0860 & 1.2424e-05 & 3.40 & 0.1143 \\
 & 128 & 9.0616e-07 & 3.52 & 0.6527 & 9.1007e-07 & 3.53 & 0.7839 \\
 & 256 & 6.2673e-08 & 3.61 & 4.6057 & 6.1772e-08 & 3.62 & 4.7521 \\
 & 512 & 3.1573e-08 & 3.26 & 34.1283 & 7.0058e-09 & 3.58 & 34.8815 \\
\end{tabular}
\end{ruledtabular}
\end{table*}

We have also tested an another unconditionally stable and fourth-order accurate Gauss-Legendre Implicit Runge–Kutta (IRK) method for the numerical computation at the initial steps. The results in the last column block of table \ref{tab:problem1_tab3} show that the scheme provides a similar accuracy as in the exact cases (see table \ref{tab:problem1_tab1} third column block) even for very small values of $\Delta t$. This behavior contrasts with the CN scheme combined with Richardson extrapolation, which did not exhibit such accuracy for small $\Delta t = 1/512$. In terms of computational efficiency, the CPU time consumed in this case is slightly higher than Richardson's extrapolation, but lower than in the case of CN with  $\delta t = \Delta t/40$. As in the earlier case of Richardson extrapolation, a fixed $9$ implicit (but not exactly CN) steps are required for the initial steps.

\begin{remark}
   It should be noted that other conditionally stable, fourth-order schemes could also be employed for the initial steps; however, in such cases, both the spatial and temporal step sizes, $\Delta x$ and $\Delta t$, must be carefully selected to satisfy the stability requirements. In the present study, since our primary focus is on the temporal order and spatial derivatives are approximated using lower-order schemes, the conventional requirement $\Delta x = \Delta t^{,p/2}$ ($p>2$, the temporal order of RBF-CN) for finding the temporal accuracy in the computation, necessitates the use of an unconditionally stable method. 
\end{remark}

\begin{table*}
\caption{\label{tab:problem1_tab4} Comparison of global error, convergence order, and CPU time of the Gaussian RBF-CN scheme with $\epsilon^2 \sim \mathcal{O}(\Delta t^3)$ when initial steps are computed using the Gauss–Legendre IRK scheme and the whole computation performed with the Gauss–Legendre IRK scheme only. }
\begin{ruledtabular}
\begin{tabular}{ccccccc}
&\multicolumn{3}{c}{RBF-CN Type III}&\multicolumn{3}{c}{Gauss-Legendre IRK}\\
$N_t$ & Global Error & Order & CPU Time &Global Error & Order & CPU Time \\ \hline
16  & 6.3065e-04 &  --  & 0.0054 & 7.8395e-05 &  --  & 0.0116 \\
32  & 3.9365e-05 & 4.00 & 0.0225 & 4.9238e-06 & 3.99 & 0.0879 \\
64  & 2.0743e-06 & 4.12 & 0.1331 & 3.0651e-07 & 4.00 & 0.7069 \\
128 & 9.8505e-08 & 4.22 & 0.8644 & 1.9211e-08 & 4.00 & 5.4154 \\
256 & 4.7746e-09 & 4.27 & 5.0852 & 1.1992e-09 & 4.00 & 46.8491 \\
512 & 7.1479e-10 & 4.06 & 37.1838 & 7.4560e-11 & 4.00 & 394.7709 \\
\end{tabular}
\end{ruledtabular}
\end{table*}


Since we observe that the Gaussian-Legendre IRK scheme performs consistently better for the initial steps calculation, irrespective of any order of approximation of $\epsilon^2$, we next employ this scheme for the entire computation and compare its performance with that of our RBF-CN scheme. In addition, the IRK method is also used to compute the initial steps in the RBF-CN scheme. Table \ref{tab:problem1_tab4} presents a comparison of the global error, convergence order, and CPU time for both schemes. The results show that both schemes achieve a similar level of accuracy; however, the RBF-CN scheme is significantly more efficient. For the smallest time step considered, $\Delta t = 1/512$, it is more than $10$ times faster and requires even less CPU time than the IRK scheme with $\Delta t = 1/256$. This demonstrates the robustness and practical efficiency of the proposed RBF-CN method, which is capable of delivering accuracy on par with a significantly less computational effort.

\begin{sidewaystable*}
\caption{\label{tab:problem1_tab5} Comparison of total number of computations for different 4th-order schemes. For the RBF-CN scheme, initial steps correspond to those used to start the scheme (first 3 $\Delta t$) and subsequent steps correspond to the main steps.}
\begin{ruledtabular}
\begin{tabular}{ccc}
Scheme & Total number of computational stages & Referred as \\ \hline
RBF-CN using IRK & $3 \times (\text{initial:}\, 2 \,\text{implicit} + 1 \,\text{explicit})$ & \\ 
& $+ 60 \times (\text{subsequent:}\, 1 \,\text{implicit} + 1 \,\text{explicit}) = 129$ & RBF-CN Type III \\
RBF-CN using CN Richardson & $3 \times (\text{initial:}\, 3 \,\text{implicit})$ & \\ 
& $+ 60 \times (\text{subsequent:}\, 1 \,\text{implicit} + 1 \,\text{explicit}) = 129$ & RBF-CN Type II \\
RBF-CN using CN with $\delta t$ & $3 \times (\text{initial:}\, 40 \,\text{implicit})$ & \\ 
& $+\, 60 \times (\text{subsequent:}\, 1 \,\text{implicit} + 1 \,\text{explicit}) = 240$ & RBF-CN Type I \\
CN Richardson & $63 \times (3 \,\text{implicit}) = 189$ & --- \\
IRK & $63 \times (2 \,\text{implicit} + 1 \,\text{explicit}) = 189$ & --- \\
RK4 & $63 \times (4 \,\text{explicit}) = 252$ & --- \\
\end{tabular}
\end{ruledtabular}
\end{sidewaystable*}


In Table \ref{tab:problem1_tab5}, we compare the total computational stages required by the RBF-CN schemes with those of other fourth-order schemes. The RBF-CN scheme with IRK requires $129$ stages in total, of which $63$ are explicit, whereas the IRK and CN Richardson schemes require $126$ and $189$ implicit stages, respectively. Although RK4 consists entirely of explicit stages and can achieve comparable computational speed, it is not unconditionally stable, unlike the RBF-CN scheme.

\begin{table*}
\caption{\label{tab:problem1_tab6} Global error, order of convergence, and CPU time of the Gaussian RBF-CN scheme with fifth order optimal $\epsilon^2$ when initial steps are computed using the Gauss–Legendre IRK scheme.}
\begin{ruledtabular}
\begin{tabular}{ccccccc}
&\multicolumn{3}{c}{RBF-CN Type III}&\multicolumn{3}{c}{Gauss-Legendre IRK}\\
&\multicolumn{3}{c}{$\epsilon^2 \sim \mathcal{O}(\Delta t^3)$}&\multicolumn{3}{c}{$\epsilon^2 \sim \mathcal{O}(\Delta t^4)$}\\
$N_t$ & Global Error & Order & CPU Time
&Global Error & Order & CPU Time \\ \hline
16  & 1.9969e-04 &  --  & 0.0116 & 2.5925e-04 &  --  & 0.0140 \\
32  & 1.6861e-05 & 3.57 & 0.0343 & 1.7972e-05 & 3.85 & 0.0553 \\
64  & 1.1408e-06 & 3.73 & 0.2088 & 3.0651e-07 & 4.86 & 0.2199 \\
128 & 5.3073e-08 & 3.95 & 1.4157 & 1.6430e-08 & 5.00 & 1.15121 \\
256 & 2.0478e-09 & 4.15 & 8.7925 & 6.7155e-10 & 4.95 & 9.1845 \\
512 & 4.6560e-10 & 3.91 & 65.8629 & 9.3123e-11 & 4.81 & 68.3532 \\
\end{tabular}
\end{ruledtabular}
\end{table*}


We next examine our schemes with the fifth-order optimal $\epsilon^2$. Table \ref{tab:problem1_tab6} compares the performance of the Gaussian RBF–CN scheme with different choices of the shape parameter scaling, namely $\epsilon^2 \sim \mathcal{O}(\Delta t^3)$ and $\epsilon^2 \sim \mathcal{O}(\Delta t^4)$. For both cases, the initial steps are computed using the Gauss–Legendre IRK method to ensure stability and accuracy at the start of the simulation. The results show that with $\epsilon^2 \sim \mathcal{O}(\Delta t^3)$, the scheme achieves an accuracy close to fourth order. In contrast, adopting the fifth-order optimal scaling $\epsilon^2 \sim \mathcal{O}(\Delta t^4)$ yields a clear improvement in the convergence behavior with the global errors decreasing at nearly fifth-order rates. The CPU times reported in the table confirm that this higher accuracy is achieved without any significant increase in computational cost. These results highlight the effectiveness of the optimal shape parameter choice in enhancing both the accuracy and efficiency of the RBF–CN scheme.

\subsection{Test Problem 2}

We next consider the diffusion equation \eqref{eq:problem1} subject to the initial and homogeneous Neumann boundary conditions
\begin{equation}
    u(x,0) = \cos(2\pi x), \qquad
    \frac{\partial u}{\partial x}(0,t) = 0, \qquad 
    \frac{\partial u}{\partial x}(1,t) = 0,
\end{equation}
so that the solution exhibits multiple zeros within the domain $x \in (0,1)$. The exact solution of this problem is given by
\begin{equation}
    u(x,t) = \cos(2\pi x) \, e^{-4\pi^2 t}.
\end{equation}

The objective of this particular test is to investigate whether the RBF-CN schemes maintain their theoretical order of accuracy in the presence of homogeneous Neumann boundary conditions and when the solution attains zeros within the domain. The time step sizes are selected such that the spatial grid resolves all the zero-crossings of the solution. On these points, the shape parameter is optimally set to zero to ensure that the RBF-CN scheme achieves at least the same order of accuracy as the standard CN. For the experiments, we employ the fourth-order optimal choice of $\epsilon^2$ for the RBF-CN schemes, while the  IRK method is used to compute the initial time steps.

\begin{figure}[htbp!]
  \centering
  (a) \hspace{3.2 in} (b) \\ 
    \includegraphics[width=0.495\textwidth]{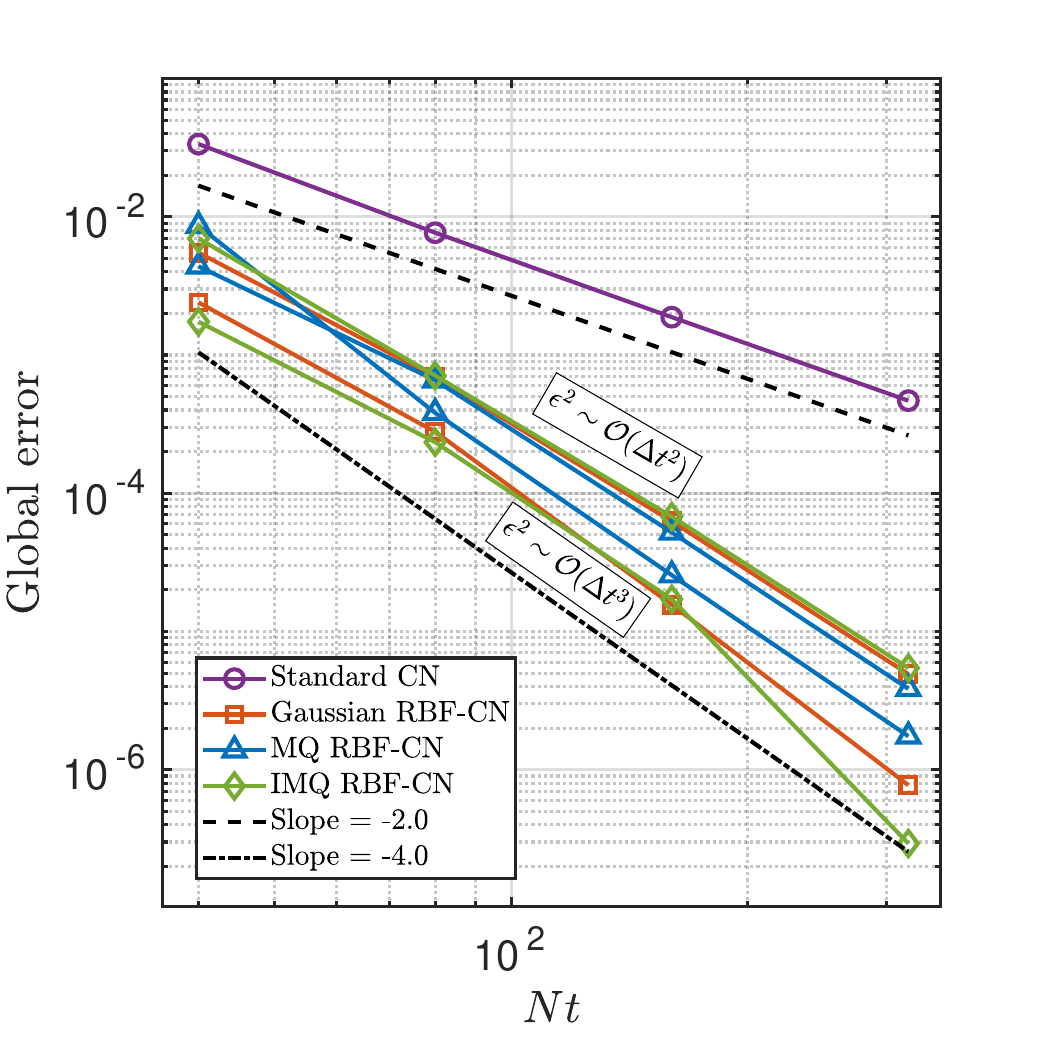}
    \includegraphics[width=0.495\textwidth]{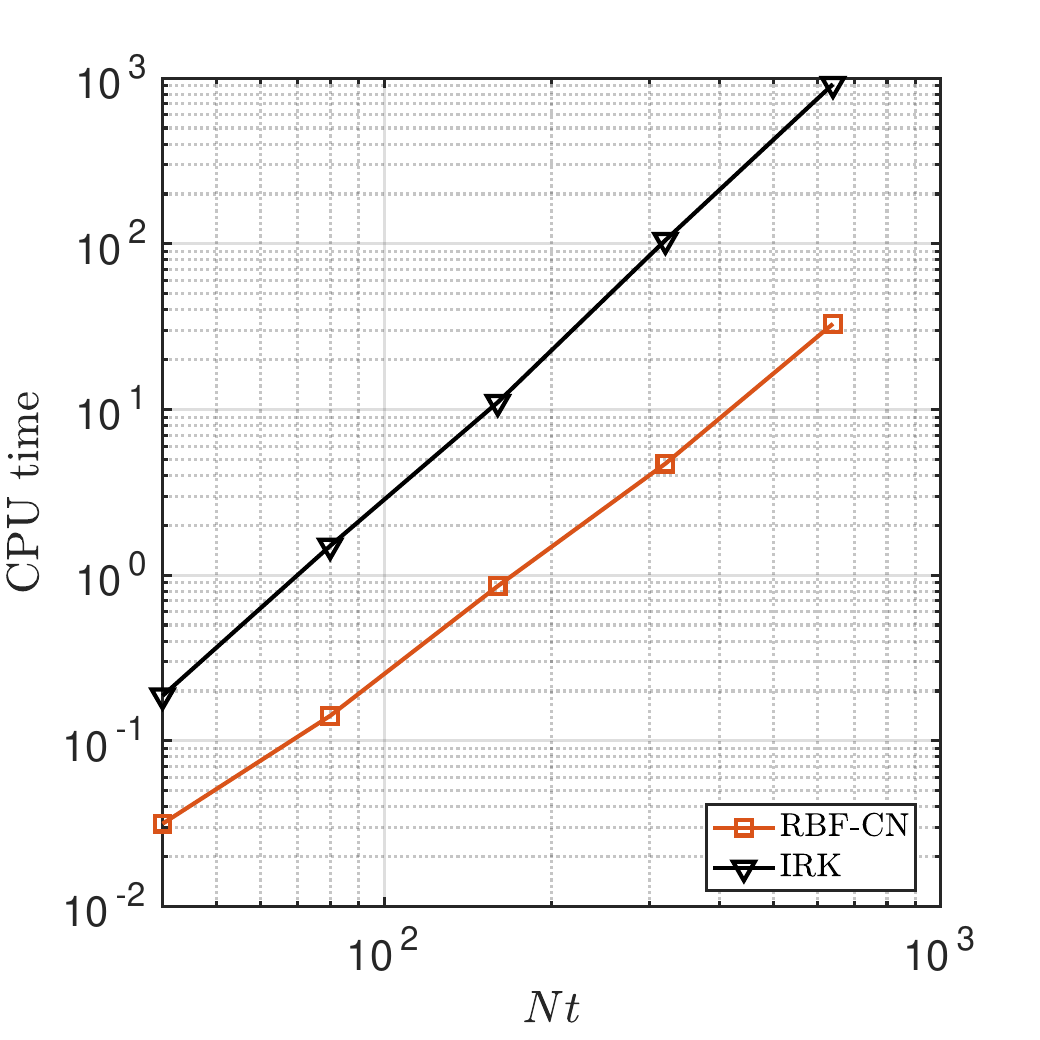}
  \caption{(a) Global convergence order of the RBF-CN schemes. (b) Comparison of the CPU time consumed by the RBF--CN scheme and the IRK scheme individually.}
  \label{fig:prob2fig1}
\end{figure}

Figure~\ref{fig:prob2fig1}(a) presents the global convergence order for three variants of the RBF-CN scheme with $\epsilon^2 \sim \mathcal{O}(\Delta t^2)$ and $\epsilon^2 \sim \mathcal{O}(\Delta t^3)$. The results reveal that, in comparison to the Dirichlet boundary condition case, the observed convergence order is slightly reduced under homogeneous Neumann boundary conditions. In particular, for the $\epsilon^2 \sim \mathcal{O}(\Delta t^2)$ approximation, the maximum order attained is approximately $3.5$, whereas for the $\epsilon^2 \sim \mathcal{O}(\Delta t^3)$ approximation, the convergence order remains close to the theoretical fourth order. Furthermore, among the three schemes examined, the MQ-based scheme achieves higher accuracy in terms of convergence order, while the Gaussian is markedly more efficient with respect to CPU time.

\begin{table*}
\caption{\label{tab:problem2_tab1}Comparison of global error, convergence order, and CPU time between the Gaussian RBF-CN scheme with $\epsilon^2 \sim \mathcal{O}(\Delta t^3)$ (using IRK for the initial steps) and the full IRK scheme.}
\begin{ruledtabular}
\begin{tabular}{ccccccc}
&\multicolumn{3}{c}{RBF-CN Type III}&\multicolumn{3}{c}{Gauss-Legendre IRK}\\
$N_t$ & Global Error & Order & CPU Time
&Global Error & Order & CPU Time \\ \hline
40  & 2.3985e-03 &  --  & 0.0316 & 5.1357e-04 &  --  & 0.1870 \\
80  & 2.7841e-04 & 3.11 & 0.1416 & 3.0768e-05 & 4.06 & 1.5015 \\
160  & 1.5625e-05 & 3.63 & 0.8591 & 1.9023e-06 & 4.04 & 11.0485 \\
320 & 7.7280e-07 & 3.90 & 4.7097 & 1.1858e-07 & 4.03 & 105.0475 \\
640 & 4.0026e-08 & 4.02 & 32.9098 & 7.4066e-09 & 4.02 & 918.7682\\
\end{tabular}
\end{ruledtabular}
\end{table*}


To further evaluate the efficiency, we compare the Gaussian RBF-CN scheme with $\epsilon^2 \sim \mathcal{O}(\Delta t^3)$, initialized using Gauss–Legendre IRK steps, against the full Gauss–Legendre IRK method. The quantitative results are reported in Table \ref{tab:problem2_tab1}, which presents the global error, convergence order, and CPU time for both schemes. The corresponding CPU time comparison is additionally depicted in Figure \ref{fig:prob2fig1}(b), providing a clearer illustration of the scaling behavior. Both methods exhibit fourth-order convergence, thereby confirming the consistency of the RBF-CN scheme. Nevertheless, the computational cost of the IRK method increases sharply with refinement, reaching nearly two orders of magnitude higher than that of the Gaussian RBF-CN scheme at the finest discretization ($N_t=640$). By contrast, the Gaussian RBF-CN scheme attains comparable accuracy with substantially reduced CPU time, demonstrating its effectiveness as a computationally efficient alternative to the full IRK method.

\section{Conclusion}\label{sec:conclusion}

In this paper, we introduced a family of modified CN schemes based on RBF interpolation in time. By considering three different infinitely smooth RBFs, namely Gaussian, MQ, and IMQ RBFs for the interpolation, we obtained three distinct RBF-CN schemes, which are found to be structure-wise simple and similar to CN. The key feature of these formulations is their ability to recover the standard CN system when the shape parameter $\epsilon^2$ involved within them is set to zero. Although the formulation does not inherently require linearity, our analysis primarily focused on the application of the schemes to linear parabolic PDEs. Consistency analysis shows that the proposed methods can achieve arbitrarily high temporal accuracy when the shape parameter $\epsilon^2$ is chosen appropriately. It is also established that the schemes always maintain at least second-order accuracy, regardless of the choice of $\epsilon^2$. Moreover, using a simple expression for the optimal $\epsilon^2$, derived from the leading term in the local truncation error, the modified scheme becomes two orders more accurate than the classical CN method, which makes it more effective and better suited for practical applications.
The stability analysis revealed that the schemes are unconditionally stable whenever the eigenvalues of the discretized spatial operator are real and non-positive, which is the same condition required for the classical CN method for the general case.

In the numerical experiment, we tested the robustness and efficiency of the proposed schemes under two different boundary conditions. In both cases, the numerical results confirmed the predicted theoretical accuracy of the proposed schemes. Since achieving higher-order accuracy with the RBF–CN schemes necessitates additional initial steps to determine the optimal value of $\epsilon^2$, we investigated several approaches to address this start-up difficulty. For the fourth-order case, CN with refined time steps, CN with Richardson extrapolation, and the Gauss–Legendre IRK method were tested. Among these, the Gauss–Legendre IRK method consistently produced highly accurate initial values across different test problems. A general principle was established that, to construct an $m$-th order accurate RBF–CN scheme, the schemes employed for the initial steps must also be at least of order $m$.

All three RBF schemes exhibited similar accuracy and required comparable CPU time, but the Gaussian RBF showed more consistent performance over the problems tested. We also examined the effect of different order of approximation of the parameter $\epsilon^2$ and found that the RBF-CN scheme is only accurate in the order of $m$ when $\epsilon^2$ is approximated by at least $\mathcal{O}(\Delta t^{m-2})$ or $\mathcal{O}(\Delta t^{m-1})$. Finally, a comparison of the RBF–CN scheme with the full Gauss–Legendre IRK method demonstrated that the proposed schemes achieve comparable convergence rates while requiring significantly lower computational cost, thereby providing a more efficient alternative for practical simulations.

Overall, the theoretical analysis and numerical evidence suggest that the proposed family of RBF–CN schemes constitutes a flexible and efficient class of time discretization methods for parabolic PDEs. The ability to achieve arbitrarily high temporal accuracy combined with unconditional stability and reduced computational overhead compared to full implicit Runge–Kutta methods underscores their potential for application in large-scale scientific computing problems. 

\begin{acknowledgments}
Authors acknowledge financial support through the Core Research Grant (CRG/2023/004156) supported by the Science and Engineering Research Board, Department of Science and Technology, Government of India. 
S.P. acknowledges financial support through the Start-Up Research Grant (SRG/2021/001269), MATRICS Grant (MTR/2022/000493) from the Science and Engineering Research Board, Department of Science and Technology, Government of India, and Start-up Research Grant (MATHSUGIITG01371SATP002), IIT Guwahati. 
\end{acknowledgments}

\nocite{*}

\bibliography{apssamp}

\end{document}